\documentclass[reqno]{amsart}
\usepackage{amssymb, graphicx, graphics, color}

\hoffset= - 0.75 in

\begin{document}
\newcommand{\hmgini}[2]{h_{#1}(#2)} % retained only for compatibility with old xfig files
\newcommand{\hmg}[3]{h_{#1}(#2,#3)} % ditto
\newcommand{\comment}[1]{}
\newcommand{\defterm}[1]{{\it #1}}
\newcommand{\ideal}[1]{\left\langle{#1}\right\rangle}
\newcommand{\coef}[2]{\left[#1\right]{#2}}      % coefficient of monomial #1 in polynomial #2
\newcommand{\tensor}{\otimes}
\newcommand{\greduce}[1]{\overline{#1}}
\newcommand{\LL}{{\mathcal L}}
\newcommand{\Zz}{{\mathbb Z}}
\newcommand{\Cc}{{\mathbb C}}
\newcommand{\Nn}{{\mathbb N}}
\newcommand{\Pp}{{\mathbb P}}
\newcommand{\Qq}{{\mathbb Q}}
\newcommand{\Sym}{{\mathfrak S}}
\newcommand{\Hilb}{\mathop{\rm Hilb}\nolimits}
\newcommand{\Poin}{\mathop{\rm Poin}\nolimits}
\newcommand{\rank}{\mathop{\rm rank}\nolimits}

\newtheorem{thm}{Theorem}[section]
\newtheorem{defn}[thm]{Definition}
\newtheorem{prop}[thm]{Proposition}
\newtheorem{lem}[thm]{Lemma}
\newtheorem{cor}[thm]{Corollary}

\theoremstyle{remark}
\newtheorem{remark}[thm]{Remark}

\title[Classification of Ding's Schubert varieties]
        {Classification of Ding's Schubert varieties:  finer rook equivalence}
%\title[Classifying Ding's partition Schubert varieties]
%        {Finer rook equivalence for Ferrers boards:  classification of
%         Ding's partition Schubert varieties}
\author[Develin, Martin, and Reiner]{Mike Develin, Jeremy L. Martin, and 
Victor Reiner}
\address{Mike Develin, American Institute of Mathematics, 360 Portage Ave., Palo Alto, CA 94306-2244, USA}
\email{develin@post.harvard.edu}
\address{Jeremy L. Martin, School of Mathematics, University of Minnesota, Minneapolis, MN 55455, USA}
\email{martin@math.umn.edu}
\address{Victor Reiner, School of Mathematics, University of Minnesota, Minneapolis, MN 55455, USA}
\email{reiner@math.umn.edu}
\date{\today}

\thanks{This work was completed while the first author was visiting the Univeristy
of Minnesota.  First author supported by the American Institute of 
Mathematics.  Second author partially supported by an
NSF Postdoctoral Fellowship.  Third author partially supported by NSF grant
DMS--0245379}


\begin{abstract}
K.~Ding studied a class of Schubert varieties $X_\lambda$ 
in type A partial
flag manifolds, indexed by
integer partitions $\lambda$ and in bijection
with dominant permutations.  He observed that the
Schubert cell structure of $X_\lambda$ is indexed by maximal rook
placements on the Ferrers board $B_\lambda$, and that the
integral cohomology groups $H^*(X_\lambda;\:\Zz)$, $H^*(X_\mu;\:\Zz)$ are
additively isomorphic exactly when the Ferrers boards $B_\lambda, B_\mu$
satisfy the combinatorial condition of \textit{rook-equivalence}.

We classify the varieties $X_\lambda$ up to isomorphism, distinguishing them
by their graded cohomology rings with integer coefficients.  The crux of our approach
is studying the nilpotence orders of linear forms in 
the cohomology ring.
\end{abstract}

\maketitle

%%%%%%%%%%%%%%%%%%%%%%%%%%%%%%%%%%%%%%%%%%%%%%%%%%%%%%%%%%%%%%%%%%%%%%%%%%%%%%%%%%%%%%%
\section{Introduction}
%%%%%%%%%%%%%%%%%%%%%%%%%%%%%%%%%%%%%%%%%%%%%%%%%%%%%%%%%%%%%%%%%%%%%%%%%%%%%%%%%%%%%%%

The goal of this paper is to classify up to isomorphism a certain class
of Schubert varieties within partial flag manifolds of type $A$.
Although this is partly motivated as a first step toward the isomorphism
classification of all Schubert varieties, we choose here to explain
instead our original motivation, stemming from rook theory in combinatorics.

A {\it board} $B$ is a subset of the squares on an $N \times N$ chessboard,
and a {\it $k$-rook placement} on $B$ is a subset of $k$ squares in $B$,
no two in a single row or column.  
Kaplansky and Riordan \cite{KaplanskyRiordan} considered
the problem of when two boards $B, B'$ are {\it rook-equivalent}, that is, when
for each $k \geq 0$, the number $R_k(B)$ of $k$-rook placements is the same
as $R_k(B')$.
%  We will use the notation $B \sim B'$ for rook-equivalence.
%%% I took this out because it was only used in two other places - JLM 3/26

Foata and Sch\"utzenberger~\cite{FoataSchutzenberger}
solved the problem for the well-behaved subclass
of {\it Ferrers boards} $B_\lambda$; these are the usual Ferrers diagrams
associated to partitions\footnote{NB: we are writing our partitions with the parts
in weakly increasing order, contrary to usual combinatorial conventions,
but more convenient in this setting.}  
\begin{equation}
\label{lambda-definition}
\lambda=(0 \leq \lambda_1 \leq \ldots \leq \lambda_n)
\end{equation}
having all squares left-justified in their row, with
$\lambda_1$ squares in the bottom row, $\lambda_2$ in the next, etc.
They showed that each rook-equivalence class of Ferrers boards has
a unique representative which is a {\it strict} partition, i.e., satisfying
$\lambda_i < \lambda_{i+1}$.
Goldman, Joichi and White~\cite{GoldmanJoichiWhite} re-proved
this result by showing that $B_\lambda$ and $B_\mu$ are rook-equivalent
if and only if the multisets of
integers $\{\lambda_i - i\}_{i=1}^n$ and $\{\mu_i - i\}_{i=1}^n$ coincide.

Garsia and Remmel~\cite{GarsiaRemmel} defined {\it $q$-rook polynomials} $R_k(B_\lambda,q)$ that
$q$-count the $k$-rook placements on $B_\lambda$ by a certain ``inversion''
statistic generalizing inversions of permutations.  They also showed
that the problem of $q$-rook equivalence is the same as that
of rook equivalence.  When $\lambda_i \geq i$ for each $i$, this can be deduced
from a product formula for $R_n(B_\lambda,q)$ that
counts placements of $n$ rooks: up to a factor of $q$ it is
\begin{equation}
\label{maximal-q-rook-polynomial}
\prod_{i=1}^n [\lambda_i - i + 1]_q
\end{equation}
where $[m]_q := \frac{q^m-1}{q-1} = 1 + q + q^2 + \cdots + q^{m-1}$.

K. Ding~\cite{Ding1,Ding2} interpreted this product as the Poincar\'e series
for a certain algebraic variety $X_\lambda$ which he called
a {\it partition variety}.  Fix a standard complete flag
of subspaces 
    $$
    0 \subset \Cc^1 \subset \cdots \Cc^{N-1} \subset \Cc^N
    $$
and define
    \begin{equation}
    \label{Ding-variety-definition}
    X_\lambda:=\{ \text{flags } 
    0 \subset V_1 \subset V_2 \subset \cdots \subset V_n \subset \Cc^N: 
    \dim_\Cc V_i = i \text{ and } V_i \subset \Cc^{\lambda_i} \}.
    \end{equation}
The set $X_\lambda$ may be endowed with the structure of a smooth complex projective
variety, and (although not stated explicitly in \cite{Ding1}) is in fact 
a smooth Schubert variety inside the partial flag manifold
$X_{N^n}$, where $N^n$ denotes the rectangular board
with $n$ rows and $N$ columns.  As we shall explain below, the
Schubert varieties
arising in this way are (in the notation of~\cite[\S 10.2]{Fulton}) those of the form
$X_w$, where $w$ is a {\it 312-avoiding} permutation.
Equivalently, the fundamental cohomology class $[X_w]$ is represented by a
Schubert polynomial indexed by a {\it dominant} or {\it 132-avoiding} permutation.
(See~\cite{Fulton} for a reference on Schubert
varieties, and~\cite{Macdonald} for a detailed treatment of
Schubert polynomials.)
Ding observed that the Schubert cell structure inherited by
$X_\lambda$ has cells indexed by $n$-rook placements on $B_\lambda$,
and with the dimension of the cell governed by Garsia and Remmel's 
inversion statistic.  Since these cells are all even-dimensional, their
(co)homology is free abelian, occurring only in even dimension, and the
Poincar\'e series of $X_\lambda $ is given by the $q$-rook polynomial formula 
\eqref{maximal-q-rook-polynomial}.
From this Ding concluded \cite{Ding2} that two partition varieties
$X_\lambda, X_\mu$ have {\it additively} isomorphic (co)homology groups
if and only if $B_\lambda$ and $B_\mu$ are rook-equivalent.

It is natural to ask when two such 
Ding partition varieties $X_\lambda, X_\mu$ 
have isomorphic (graded) cohomology {\it rings}, or even
when they are isomorphic {\it as varieties}.
The main result of this paper is that the answers to both questions are the same.
We make use of recent results of Gasharov and 
the third author \cite{GasharovReiner},
giving simple explicit cohomology ring 
presentations\footnote{It is amusing that these cohomology ring presentations for
Schubert varieties are often derived for the purposes of enumerative
geometry (Schubert calculus), but are used here for a different  
classical topological purpose, namely distinguishing non-homeomorphic spaces.}
for a more general class
of Schubert varieties in partial flag manifolds (those defined by
a conjunction of inclusion conditions of the forms 
$\Cc^j \subset V_i$ and $V_i\subset \Cc^j$).

To state our main result, we first note one trivial source of
isomorphisms among the partition varieties $X_\lambda$.
We assume throughout that $\lambda_i \geq i$ for every $i$, for otherwise
$X_\lambda=\emptyset$.  However, if $\lambda_k=k$ for some $k$,
then the condition $V_k \subset \Cc^k$ with $\dim_\Cc V_k=k$ forces $V_k = \Cc^k$,
so that $X_\lambda$ is isomorphic to $X_{\lambda^{(1)}}
\times X_{\lambda^{(2)}}$, where
$$
\begin{aligned}
\lambda^{(1)} &= (\lambda_1,\,\ldots,\,\lambda_{k-1}),\\
\lambda^{(2)} &= (\lambda_{k+1}-k,\,\ldots,\,\lambda_n-k).
\end{aligned}
$$
Here if $k=n$, so that $\lambda_n=n$, there is no
partition $\lambda^{(2)}$ and we simply note that $X_\lambda \cong X_{\lambda^{(1)}}$.

\begin{figure}
\begin{center}
\resizebox{5.4cm}{4.8cm}{\includegraphics{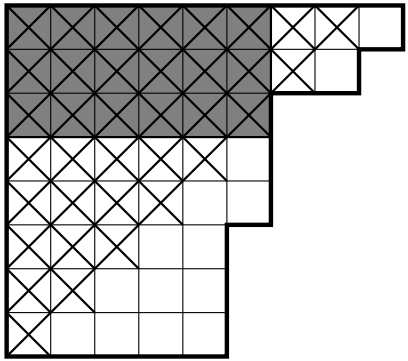}} % revised JLM 2/25
\end{center}
\caption{\label{lambda-example}
A decomposable partition $\lambda$.  The unshaded regions are $\lambda^{(1)}$
and $\lambda^{(2)}$.}
\end{figure}

%\noindent
Say that $\lambda$ is {\it decomposable} if this occurs (i.e., if $\lambda_k=k$
for some $k$), and {\it indecomposable} otherwise.
For example, the partition $\lambda=(5,5,5,6,6,6,8,9)$ shown in Figure~\ref{lambda-example}
is decomposable since $\lambda_6=6$.  In this case, one has
$\lambda^{(1)} = (5,5,5,6,6)$ and $\lambda^{(2)} = (2,3)$, as shown in the figure.

Iterating this, one can decompose $\lambda$ 
into a multiset of indecomposable partitions $\{\lambda^{(i)}\}_{i=1}^r$, 
which we will call its \defterm{indecomposable components}, such that 
    \begin{equation} \label{decompose-variety}
    X_\lambda \cong X_{\lambda^{(1)}} \times \cdots \times X_{\lambda^{(r)}}.
    \end{equation}

%\begin{figure}
%\input{gen-lambda.pstex_t}
%\caption{\label{gen-lambda}
%A partition $\lambda$. This decomposes into the partitions $\lambda_1$ and 
%$\lambda_2$, effecting a tensor product decomposition of the corresponding 
%cohomology rings.}
%\end{figure}

Our main result is that the Schubert varieties $X_\lambda$ are determined
up to isomorphism by these multisets of indecomposable components.  
It should be compared with the result of Goldman, Joichi and White
\cite{GoldmanJoichiWhite}, which can now be rephrased:  
the varieties $X_\lambda$ are 
determined up to additive (co-)homology isomorphism by the multisets
of numbers $\{\lambda_i-i\}$.

\begin{thm} \label{main-theorem}
The following are equivalent for two partitions,
$\lambda=(\lambda_1,\ldots,\lambda_m)$ and $\mu=(\mu_1,\ldots,\mu_{m'})$:
\begin{enumerate}
\item[(i)] The multisets of indecomposable components, 
$\{\lambda^{(i)}\}_{i=1}^r$ and $\{\mu^{(i)}\}_{i=1}^{r'}$,
are identical.
\item[(ii)] There is an isomorphism $X_\lambda \cong  X_\mu$ of varieties.
\item[(iii)] There is a graded isomorphism of integer cohomology rings
$H^*(X_\lambda;\:\Zz) \cong H^*(X_\mu;\:\Zz)$.
\end{enumerate}
\end{thm}

The implications (i) $\implies$ (ii) $\implies$ (iii) are clear; the
hard part is to show that (iii) implies (i).  It turns out that
the key to this implication lies in understanding the \defterm{nilpotence orders}
of cohomology elements $f \in H^2(X_\lambda;\:\Zz)$; that is,
the least $k$ for which $f^k=0$.  

In Section~\ref{presentation-review}, we review some of Ding's results,
and re-prove somewhat more directly the presentation for 
$H^*(X_\lambda;\:\Zz)$ from \cite{GasharovReiner}.
The three sections that follow are the technical heart of the paper,
categorizing elements in $H^2(X^\lambda;\:\Zz)$ of minimal
nilpotence order.  We begin in Section~\ref{GB-section}
by setting up some Gr\"obner basis
machinery that we shall use throughout (for a general reference on
Gr\"obner basis theory, see~\cite{CoxLittleOShea}).  Section~\ref{G/B-section} deals
with nilpotents in the cohomology of the complete flag variety
(that is, when $\lambda$ is a square Ferrers board) and
Section~\ref{general-nilpotence-section} treats the case of arbitrary $X_\lambda$.
Using these tools, we prove in Section~\ref{indecomposable} that
an indecomposable partition $\lambda$ may be recovered from the structure of
$H^*(X_\lambda;\:\Zz)$ as a graded $\Zz$-algebra.  Finally, in
Section~\ref{decomposable}, we show that in the general case, $H^*(X_\lambda;\:\Zz)$
has an essentially unique decomposition as a tensor product of graded
$\Zz$-algebras, whose factors correspond to the indecomposable components of
the partition $\lambda$.

It is curious that this unique tensor decomposition fails if instead of
the integer cohomology ring $H^*(X_\lambda;\:\Zz)$ one takes cohomology
with coefficients in a ring where $2$ is invertible; 
see Remark~\ref{invertible-2} below.

%%%%%%%%%%%%%%%%%%%%%%%%%%%%%%%%%%%%%%%%%%%%%%%%%%%%%%%%%%%%%%%%%%%%%%%%%%%%%%%%%%%%%%%
\section{Review of $X_\lambda$ and the presentation of $H^*(X_\lambda;\:\Zz)$}
\label{presentation-review}
%%%%%%%%%%%%%%%%%%%%%%%%%%%%%%%%%%%%%%%%%%%%%%%%%%%%%%%%%%%%%%%%%%%%%%%%%%%%%%%%%%%%%%%

For the sake of completeness, and also to collect facts for future
use, we begin by re-proving some of Ding's results from \cite{Ding1}, and
re-derive somewhat more directly the
presentation given in~\cite{GasharovReiner} for the cohomology ring of $X_\lambda$.
Throughout this paper, all cohomology groups and rings are taken with integer
coefficients unless otherwise specified.
We begin by identifying the Schubert varieties that arise as
Ding's varieties $X_\lambda$.  (See~\cite[\S 10.6]{Fulton}
for more information on Schubert varieties, and~\cite{Macdonald} for a
detailed treatment of Schubert polynomials.)

Let $\Sym_N$ be the symmetric group
of permutations of $\{1,\dots,N\}$, and let $\Sym_{\{n+1,n+2,\ldots,N\}}$ be
the subgroup of permutations fixing $\{1,\ldots,n\}$ pointwise.
Consider the partial flag variety
    $$
    X_{N^n} ~=~ \big\{\text{flags }0 \subset V_1 \subset \cdots \subset V_n \subset \Cc^N:\ \ 
    \dim V_i = i\big\}.
    $$
Let $w=w_1\dots w_n \in \Sym_N$ be a permutation which is
a maximum-length representative
for its coset in $\Sym_N/\Sym_{\{n+1,n+2,\ldots,N\}}$.  The
corresponding Schubert variety $X_w \subset X_{N^n}$ is defined to be
    $$
    X_w ~=~ \big\{\text{flags }0 \subset V_1 \subset \cdots \subset V_n \subset \Cc^N:\ \ 
      \dim V_i = i,\ \dim V_i \cap \Cc^j \geq \# \{k \leq i: w_k \leq j\} \big\}.
    $$

Let $\lambda$ be a partition of the form~\eqref{lambda-definition}, and let
$N = \lambda_n$.  It is easy to check that Ding's variety $X_\lambda$
coincides with the Schubert variety $X_w \subset X_{N^n}$, where $w$ is the
unique permutation given by the recursive rule
    $$w_i ~=~ \max\left( \{1,\,\dots,\,\lambda_i\} \setminus
      \{w_1,\,\dots,\,w_{i-1}\} \right).$$
Note that if $n=N$, then $w$ corresponds to the
maximal rook placement on the Ferrers board $B_\lambda$ given by the
following algorithm: let $i$ increase from $1$ to $n$, and for each $i$,
place a rook in row $i$ and column $w_i$, where $w_i$ is the
rightmost square in row $i$ whose column does not already contain a rook.
For instance, if $\lambda = (2,4,4,5,5)$, then $w=24351 \in \Sym_5$.
(If $n<N$, then we must first augment $\lambda$ with $N-n$ additional rows
of length $\lambda_n$.)
It is not hard to verify that the permutations $w$ obtained in this way are
exactly those which are {\it 312-avoiding}; that is, there do not exist
$i,j,k$ for which $i<j<k$ and $w(i)>w(k)>w(j)$.  Equivalently, the cohomology class
$[X_w] \in H^*(X_{N^n})$ is represented by a
Schubert
polynomial which is a single monomial, namely the Schubert polynomial indexed by the
dominant (or 132-avoiding) permutation $w_0w$, where $w_0$ is the unique permutation
of maximal length.
(We thank Ezra Miller for discussions
clarifying these points.)

Because $X_\lambda$ is a Schubert variety, it comes equipped with
a Schubert cell decomposition, having cells in only even real dimensions.
As observed by Ding, this has important consequences:

\begin{thm}[Ding~\cite{Ding1}]
\label{Ding's-thm}
The integral cohomology ring $H^*(X_\lambda;\:\Zz)$ is free abelian (that is, it has
no torsion), is nonzero only in even homological degrees, and has
Poincar\'e series
$$
\Poin(X_\lambda, q) := \sum_{i \geq 0} q^i \rank_\Zz H^{2i}(X_\lambda;\:\Zz)
= \prod_{i=1}^n [\lambda_i - i + 1]_q.
$$
\end{thm}

\begin{proof}
The cohomology is free abelian and concentrated in even degrees
because the Schubert cell decomposition for the Schubert
variety $X_\lambda$ has cells only in even dimensions.  

For the assertion about the Poincar\'e series, we will induct on $n$.
The map 
$$
\begin{matrix}
X_\lambda & \rightarrow & \Pp(\Cc^{\lambda_1}) \cong \Pp_\Cc^{\lambda_1-1} \\
\{V_i\}_{i=1}^n & \mapsto & V_1
\end{matrix}
$$
is an (algebraic) fiber bundle, with fiber isomorphic to $X_\nu$,
where 
$$
\nu=(\nu_1,\ldots,\nu_{n-1}) = (\lambda_2-1,\ldots,\lambda_n-1)
$$
is the partition obtained by removing the first row and column from $\lambda$
(see Figure~\ref{nu-example}).
\begin{figure}
\begin{center}
\resizebox{6.0cm}{3.6cm}{\includegraphics{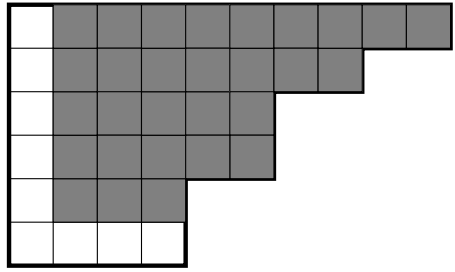}}
\end{center}
\caption{\label{nu-example}
A partition $\lambda$ and the subpartition $\nu$ (shaded)
such that $H^*(X_\lambda)\:/\ideal{x_1} = H^*(X_\nu)$.}
\end{figure}
The Leray-Serre spectral sequence is particularly simple in this situation,
because both base and fiber are simply-connected (again due to the
Schubert cell decomposition) and have homology concentrated in even
dimension.  This causes the spectral sequence to
degenerate at the $E^1$-page, yielding 
$$
\Poin( X_\lambda, q ) = \Poin( X_\nu, q ) \cdot \Poin( \Pp_\Cc^{\lambda_1-1}, q ).
$$
The assertion about $\Poin( X_\lambda, q )$ now follows by induction on $n$,
using the fact that $[m]_q = \Poin( \Pp^{m-1}_\Cc, q)$.
\end{proof}

We now set about deriving the presentation for $H^*(X_\lambda)$.
To this end, we recall
Borel's picture for the cohomology of the complete flag manifold
$GL_N(\Cc)/B = X_{N^N}$ and the partial flag manifold $X_{N^n}$;
see \cite[Chapter 10, \S3, \S6]{Fulton}.
We will use the following notation for symmetric functions in various sets of
variables.
For integers $1\leq i \leq j \leq N$ and $m \geq 0$, define
the $m^{th}$ \defterm{elementary} and \defterm{complete homogeneous}
symmetric functions, respectively, by
\begin{align*}
e_m(i,j) &\;:=\; e_m(x_i,x_{i+1},\ldots,x_j), \\
h_m(i,j) &\;:=\; h_m(x_i,x_{i+1},\ldots,x_j), \\
\intertext{where}
e_m(N) &\;:=\; e_m(1,N) \;=\; e_m(x_1,\ldots,x_N) ~=~
   \sum_{1 \leq i_1 < \cdots < i_m \leq N} x_{i_1} \cdots x_{i_m}, \\
h_m(N) &\;:=\; h_m(1,N) \;=\; h_m(x_1,\ldots,x_N) ~=~
   \sum_{1 \leq i_1 \leq \cdots \leq i_m \leq N} x_{i_1} \cdots x_{i_m} .
\end{align*}

According to Borel's picture, $H^*(X_{N^N}) \cong \Zz[x_1,\ldots,x_N]/J$,
where $J=\ideal{e_1(N),\ldots,e_N(N)}$ is the ideal generated by all symmetric functions
of positive degree, and where $x_i$ represents the negative of $c_1(\LL_i)$,
the first Chern class of the line bundle $\LL_i$ on $GL_N(\Cc)/B$ whose fiber
over the flag $\{V_i\}_{i=0}^N$ is $V_i/V_{i-1}$.

Furthermore, the surjection $X_{N^N} \rightarrow X_{N^n}$ which forgets
the subspaces of dimension greater than $n$ in a complete flag induces
a map $H^*(X_{N^n}) \rightarrow H^*(X_{N^N})$ which turns out to be injective,
and the image of $H^*(X_{N^n})$ is identified with the invariant subring
$H^*(X_{N^N})^{\Sym_{\{n+1,n+2,\ldots,N\}}}$.  This invariant subring
may be presented as $S/J'$, where
$$
S ~=~ \Zz[x_1,\ldots,x_N]^{\Sym_{\{n+1,n+2,\ldots,N\}}}
~=~ \Zz[x_1,\ldots,x_n,\; e_1(n+1,N),\ldots,e_{N-n}(n+1,N)]
$$
and $J'=\ideal{e_1(N),\ldots,e_N(N)}$ is the ideal of $S$
with the same generators as $J$.

The relations in the ideals $J$ and $J'$ induce further relations
among various symmetric functions, which we record here for future
use.

\begin{prop}[cf.~{\cite[p.~163,~eqn.~(4)]{Fulton}}]
\label{Fulton-trick}
For every $m \in \{1,2,\dots,N\}$ and $j \geq 0$, one has
$$
h_j(m) \equiv (-1)^j e_j(m+1,N) \pmod{J}.
$$
\end{prop}
\begin{proof}
$$
\prod_{i=1}^m (1+x_i t) \prod_{i=m+1}^N (1+x_i t) =
\prod_{i=1}^N (1+x_i t) = \sum_{j=0}^n e_j(N) t^j
\equiv 1 \pmod{J}.
$$
Hence
$$
\sum_{j=0}^{\infty} (-1)^j h_j(m) t^j
= \prod_{i=1}^m (1+x_i t)^{-1}
\equiv \prod_{i=m+1}^N (1+x_i t)
= \sum_{j=0}^{N-m} e_j(m+1,N) t^j \pmod{J}.
$$
Now comparing coefficients of powers of $t^j$ yields the desired equality.
\end{proof}

We now give the general presentation for the integral cohomology of $X_\lambda$
(as pointed out in~\cite[Remark 3.3]{GasharovReiner}).

\begin{thm}
\label{cohomology-presentation}
Let $\lambda$ be a partition with 
$1 \leq \lambda_1 \leq \cdots \leq \lambda_n = N$ and $\lambda_i \geq i$ for all $i$.
Let 
$$
R^\lambda:=\Zz[x_1,\ldots,x_n]/I_\lambda
$$
where $I_\lambda := \ideal{h_{\lambda_i-i+1}(i) : 1 \leq i \leq n}$.

Then there is a (grade-doubling) ring isomorphism
$$
R^\lambda \rightarrow H^*(X_\lambda;\:\Zz)
$$
sending $x_i$ 
to $-c_1(\LL_i)$.  Here $\LL_i$ is the same line bundle as above,
but restricted to $X_\lambda$ from the partial flag manifold $X^{N^n}$.
\end{thm}

\begin{proof}
The obvious inclusion $X_\lambda \hookrightarrow X_{N^n}$
induces a map $H^*(X_{N^n}) \rightarrow H^*(X_\lambda)$.
This ring map is surjective, because $X_\lambda$ inherits from $X_{N^n}$
a decomposition into Schubert cells, and
the dual cocycles to these (even-dimensional) cells additively 
generate the cohomology in each case.

There are further relations on the Chern classes $x_i$ in $H^*(X_\lambda)$
due to the conditions $V_i \subset \Cc^{\lambda_i}$.  Specifically,
the bundle $\Cc^N / V_i$ on $X_\lambda$ will have the
same Chern classes as the direct sum 
$\Cc^N/\Cc^{\lambda_i} \oplus \Cc^{\lambda_i}/V_i$, in which 
$\Cc^N/\Cc^{\lambda_i}$ is a trivial bundle.  Thus when 
restricted to $X_\lambda$, the bundle $\Cc^N / V_i$ will have the
same Chern classes as the bundle $\Cc^{\lambda_i}/V_i$ of rank
$\lambda_i-i$.  Hence its Chern classes $c_m = \pm e_m(i+1,N)$ 
for $m > \lambda_i-i$ inside $H^*(X_\lambda)$ must vanish.
Consequently, we have a surjection of rings
    \begin{equation}
    \label{e-presentation}
     \Zz[x_1,\ldots,x_n,e_1(n+1,N),e_2(n+1,N),\ldots,e_{N-n}(n+1,N)]/J_\lambda
    ~\rightarrow~ H^*(X_\lambda)
    \end{equation}
where 
    $$
    J_\lambda := J' + 
    \ideal{ e_j(i+1,N) ~:~ 1 \leq i \leq n\text{ and } j > \lambda_i-i }.
    $$

We now manipulate the quotient ring 
$\Zz[x_1,\ldots,x_N]^{\Sym_{\{n+1,n+2,\ldots,N\}}}/J_\lambda$
on the left of \eqref{e-presentation}.  We use 
Proposition~\ref{Fulton-trick} to draw two conclusions:
\begin{enumerate}
\item[(i)] Applying Proposition~\ref{Fulton-trick} with $m=n$ shows that
$H^*(X_{N^n})$ and $H^*(X_\lambda)$ are generated as algebras
by $x_1,\ldots,x_n$, since their generators of the form
$e_i(n+1,N)$ can be expressed modulo $J'$ as (symmetric) polynomials
in $x_1,\ldots,x_n$.
\item[(ii)] Applying it with $m=i$ for $1 \leq i \leq n$ shows
that $h_{\lambda_i-i+1}(i)=0$ in $H^*(X_\lambda)$,
because for each $j \geq \lambda_i-i$, $h_j(i)$ is congruent modulo
$J'$ to $\pm e_j(i+1,N)$.
\end{enumerate}

Consequently, there is a surjection of rings
\begin{equation}
\label{presentation-as-quotient}
\Zz[x_1,\ldots,x_n] / \ideal{ h_{\lambda_i-i+1}(i): 1 \leq i \leq n }
~\rightarrow~ H^*(X_\lambda).
\end{equation}

On the other hand, the set
    $$\{h_{\lambda_i-i+1}(i): 1 \leq i \leq n\}$$
is a Gr\"obner basis for $I_\lambda$
with respect to the lexicographic term order on $\Zz[x_1,\ldots,x_n]$
given by $x_1 < \cdots < x_n$.
Indeed, the initial term of $h_{\lambda_i-i+1}(i)$ is $x_i^{\lambda_i-i+1}$,
so these generators have pairwise relatively
prime, monic initial terms. 
Consequently, the quotient ring on the
left of \eqref{presentation-as-quotient} is a free $\Zz$-module
of rank $\prod_{i=1}^n (\lambda_i - i + 1)$, with $\Zz$-basis
given by the {\it standard monomials} (those divisible by none
of the initial terms), namely 
$\{ x_1^{a_1}\cdots x_n^{a_n}: a_i \leq \lambda_i-i \}$.  
Since Theorem~\ref{Ding's-thm} implies
that $H^*(X_\lambda)$ is a free $\Zz$-module of the same rank,
the surjection \eqref{presentation-as-quotient} must be an isomorphism.
\end{proof}

For example, if $\lambda$ is the partition shown in Figure~\ref{lambda-example},
then the Gr\"obner basis for $I_\lambda$ is
    $$h_5(1), \ h_4(2), \ h_3(3), \ h_3(4), \ h_2(5), \ h_1(6), \ h_2(7), \ h_2(8).$$

%Note that the degree of the $i^{th}$ generator,
%namely $\lambda_i-i+1$, equals the number of unmarked boxes in the $i^{th}$ row
%of Figure~\ref{lambda-example}.

The previous proof shows that $I_\lambda$ is the {\it elimination
ideal}
$$
I_\lambda = \Zz[x_1,\ldots,x_n] \cap J_\lambda.
$$
This observation has some useful corollaries, which
can also be proved by direct combinatorial/algebraic arguments avoiding
any use of geometry.
%%
%% I removed the next two corollaries because I couldn't find explicit
%% references to them later ...  Vic 2/18/04
%\begin{cor} \label{square}
%Let $n \in \Nn$.  Then
%\begin{eqnarray*}
%\ideal{e_1(n), e_2(n), \ldots, e_n(n)}
%&=& \ideal{h_1(n), h_2(n), \ldots, h_n(n)}\\
%&=& \ideal{h_1(n), h_2(n-1), \ldots, h_n(1)},
%\end{eqnarray*}
%and the cohomology of the complete flag variety
%$G/B = X^{N^N}$ may be presented as the quotient
%of $\Zz[x_1,\dots,x_n]$ by any of these ideals.
%\end{cor}
%
%\begin{cor} (* \label{rectangle} *)
%Let $n,a \in \Nn$.  Then
%\begin{eqnarray*}
%\ideal{h_{a+1}(n), h_{a+2}(n), \dots, h_{a+n}(n)}
%&=& \ideal{h_{a+1}(n), h_{a+2}(n-1), \dots, h_{a+n-1}(2), h_{a+n}(1)},
%\end{eqnarray*}
%and the cohomology of the partial flag variety
%$X^{N^n}$, where $N=n+a$, may be presented as the quotient
%of $\Zz[x_1,\dots,x_n]$ by either of these ideals.
%\end{cor}
%
The first corollary is the algebraic manifestation of the
(surjective) map $R^\lambda \to R^\mu$ induced by the inclusion of 
Schubert varieties $X_\lambda \to X_\mu$.

\begin{cor} \label{contain-epi}
Let $\lambda$ and $\mu$ be partitions, both with $n$ nonzero rows, 
such that $\lambda \supset \mu$.  

Then $I_\lambda \subset I_\mu$, and consequently, $R^\lambda$ is
a quotient of $R^\mu$.
\end{cor}

\begin{proof}
By definition of $J_\lambda$, one has
$J_\lambda \subset J_\mu$ in this situation.
\end{proof}

\begin{cor}
\label{variable-symmetry}
If $\lambda_i=\lambda_{i+1} = \cdots = \lambda_j$ for some $i < j$
then the ideal $I_\lambda$ is invariant under
permutations of the variables $x_i,x_{i+1},\ldots,x_j$.
\end{cor}

\begin{proof}
It suffices to show that $J_\lambda$ has this same invariance.
Note that the generators for $J_\lambda$ of the form
$$
e_m(i'+1,N) \quad \text{for } i \leq i' < j \text{ and } m > \lambda_{i'}-i'
$$
are all redundant, as they lie in the ideal generated
by $\{e_m(j+1,N): m > \lambda_j-j\}$.  The latter generators,
and all other generators of $J_\lambda$, are symmetric in $x_i,x_{i+1},\ldots,x_j$.
\end{proof}

%%%%%%%%%%%%%%%%%%%%%%%%%%%%%%%%%%%%%%%%%%%%%%%%%%%%%%%%%%%%%%%%%%%%%%%%%%%%%
\section{Two reduced Gr\"obner bases}
\label{GB-section}
%%%%%%%%%%%%%%%%%%%%%%%%%%%%%%%%%%%%%%%%%%%%%%%%%%%%%%%%%%%%%%%%%%%%%%%%%%%%%

This section examines the Gr\"obner bases for $I_\lambda$
for two extreme cases of indecomposable partitions.  In both cases,
one can describe the (unique) reduced Gr\"obner basis, which will be used
in an essential way later in the paper.  We assume some
familiarity with ``Gr\"obner basics'' on the reader's part;
a good reference for this topic is~\cite{CoxLittleOShea}.

We begin with some notation regarding Gr\"obner reduction.
Since the generators $\{ h_{\lambda_i-i+1}(i) : 1 \leq i \leq n \}$
form a Gr\"obner basis for $I_\lambda$ with respect to a
lexicographic monomial ordering in which $x_1 < \ldots x_n$, we can
compute in the quotient $R^\lambda$ by reducing polynomials
modulo this Gr\"obner basis.  For a polynomial $f \in \Zz[x_1,\ldots,x_n]$,
we will denote by $\greduce{f}$ this \defterm{standard form} of $f$.
That is, $\greduce{f}$ is the unique $\Zz$-linear combination
of standard monomials $\{x_1^{a_1}\cdots x_n^{a_n}: a_i \leq \lambda_i-i\}$
which is congruent to $f$ modulo $I_\lambda$.  Given a standard monomial $M$,
we denote by $[M]\greduce{f}$ the coefficient of $M$ in $\greduce{f}$.
(This is well-defined, because the standard monomials form a basis
for $\Zz[x_1,\dots,x_n]/I_\lambda$ as a free $\Zz$-module.)

Let $\lambda=(\lambda_1 \leq \cdots \leq \lambda_n)$ and
for some fixed $m \leq n$, let $\mu=(\lambda_1,\ldots,\lambda_m)$.
Then the fact that we are using a lexicographic order
to perform reductions has the following easy consequence
(see also \cite[\S 3.1]{CoxLittleOShea}), which will be used frequently.
It can be viewed as an algebraic consequence of the
fibration $X_\lambda \rightarrow X_\mu$ that forgets the subspaces of 
dimension greater than $m$ in a flag, which happens to induce an injective
map $H^*(X_\mu) \to H^*(X_\lambda)$.

\begin{prop}
\label{Leray-Hirsch}
Let $\lambda$ and $\mu$ be related as above.  Suppose that
$f$ in $\Zz[x_1,\ldots,x_n]$ lies in
some subalgebra $\Zz[x_1,\ldots,x_m]$, where $m \leq n$.

Then the images of $f$ in $R^\lambda$ and $R^\mu$ have the
same standard form $\greduce{f}$.
\end{prop}

Our first extreme case arises when $\lambda$ is an indecomposable
partition with $\lambda_i=p$, and $\mu \subset \lambda$
is the smallest indecomposable partition
having $\mu_i=p$, namely $\mu=(2,3,\ldots,i-1,i,p)$.

\begin{prop}
\label{extreme-with-given-part-GB}
Let $\mu=(2,3,\ldots,i-1,i,p)$. With respect
to lexicographic order on $\Zz[x_1,\ldots,x_m]$ with $x_1 < \dots < x_m$,
the ideal $I_\mu$ has reduced Gr\"obner basis
    \begin{equation}
    \label{putative-GB}
    \{x_1h_1(1),\; x_2h_1(2),\; \ldots,\; x_{i-1}h_1(i-1),\;
      x_i^{p-i+1}+x_i^{p-i} h_1(i-1)\}.
    \end{equation}
\end{prop}

\begin{proof}
It is easy to see that the elements of \eqref{putative-GB} form a reduced
Gr\"obner basis with respect to the lexicographic order for whatever ideal
they generate.  We observe that this ideal may also be presented as
    $$\ideal{h_2(1),\; h_2(2),\; \ldots,\; h_2(i-1),\; x_i^{p-i+1}+x_i^{p-i}h_1(i-1)}.$$
We will show that this ideal is exactly $I_\mu$.
By Theorem~\ref{cohomology-presentation},
    $$I_\mu ~=~ \ideal{h_2(1),h_2(2),\ldots,h_2(i-1),h_{p-i+1}(i)},$$
so it remains only to show that
$h_{p-i+1}(i)$ and $x_i^{p-i+1}+x_i^{p-i}h_1(i-1)$ are congruent
modulo the ideal $\ideal{h_2(1),h_2(2),\ldots,h_2(i-1)}$.
Since $h_{p-i+1}(i) = \sum_{j=1}^{p-i+1} x_i^j h_{p-i-j+1}(i-1)$,
this congruence is immediate from the fact that
$$
h_m(\ell) \in \ideal{h_2(1), h_2(2), \ldots, h_2(\ell)} \text{ for }m \geq 2,
$$
which is easily proven by double induction on $m$ and $\ell$
via the identity $h_m(\ell) = x_\ell h_{m-1}(\ell) + h_m(\ell-1)$.
\end{proof}

\begin{figure}
\begin{center}
\resizebox{6.0cm}{3.6cm}{\includegraphics{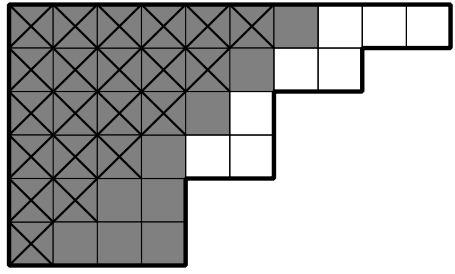}}
\end{center}
\caption{\label{core-figure}
An indecomposable partition $\lambda$ and its core subpartition $\mu$ (shaded).}
\end{figure}

Our second extreme case arises when 
$\lambda$ is an indecomposable partition with $n$ rows.  Let $k=\lambda_1$, and
let $\mu$ be the smallest indecomposable partition with $n$ rows and $\mu_1=k$.
That is,
    \begin{equation} \label{core}
    \begin{aligned}
    & \mu_1 = \mu_2 = \dots = \mu_{k-1} = k, \\
    & \mu_i = i+1 \quad \text{for } k \leq i \leq n.
    \end{aligned}
    \end{equation}
\noindent
Then $\mu$ is a subpartition\footnote{
For the purposes of this paper, the statement
``$\mu$ is a subpartition of $\lambda$''
means that $\mu_i\leq\lambda_i$ for all rows $\mu_i$ of $\mu$.
Equivalently, the Ferrers diagram of $\mu$ is contained inside
that of $\lambda$, when both are left- and bottom-justified.}
of $\lambda$, which we will call the \defterm{core} of $\lambda$.
For example, the core of $\lambda=(4,4,6,6,8,10)$ is the 
partition $\mu=(4,4,4,5,6,7)$
(see Figure~\ref{core-figure}). 

%The benefit of Corollary~\ref{contain-epi} is that we can frequently
%answer questions about $R^\lambda$ by passing to $R^\mu$.

\begin{prop} \label{grobbas}
For $k<n$, let $\lambda$ be a partition which is its own core.

Then the polynomials
    \begin{equation}\label{RGB}
    \begin{aligned}
    G_1 &= h_k(1), & G_2 &= h_{k-1}(2), &&\dots, & G_{k-1} &= h_2(k-1),\\
    G_k &= x_k h_1(k), & G_{k+1} &= x_{k+1} h_1(k+1), &&\dots, & G_n &= x_n h_1(n)
    \end{aligned}
    \end{equation}
form a reduced Gr\"{o}bner basis for $I_\lambda$ under the reverse lexicographic
term order given by $x_1<x_2<\ldots<x_n$.
\end{prop}

\begin{proof}
The initial terms of the $G_i$'s are (in order)
$x_1^k$, $x_2^{k-1}$, \ldots, $x_{k-1}^2$, $x_k^2$, \ldots, $x_n^2$.
It is evident that no initial term divides any term of any other
$G_i$.  Therefore, they are a reduced Gr\"obner basis for the ideal that they generate.

We claim that for every $r \in \{k,k+1,\dots,n\}$,
    $$\ideal{G_1,\dots,G_r} ~=~
      \ideal{h_k(1),\, h_{k-1}(2),\, \dots,\, h_2(k-1),\, h_2(k),\, \dots,\, h_2(r)}.$$
The claim is trivial for $r=k$.  For $r>k$, it follows from
induction and the observation that $h_2(r)-h_2(r-1) = x_r h_1(r) = G_r$.
In particular, the equality for $r=n$ gives $\ideal{G_1,\dots,G_r} = I_\lambda$.
\end{proof}

The form of this reduced Gr\"obner basis has the following consequence,
which we will exploit later.
\begin{cor} \label{sticky}
(``Stickiness'')  Let $\lambda$ be an indecomposable partition which is its own core,
and $k:=\lambda_1$.  Let $M$ be a monomial in $x_1,\dots,x_n$.

Then:
\begin{enumerate}
\item If $k \leq i \leq n$ and $M$ is divisible by $x_i$, then so is $\greduce{M}$.
\item If $M$ is not divisible by any of the variables $x_k,\dots,x_n$,
then neither is $\greduce{M}$.
\end{enumerate}
\end{cor}

\begin{proof}
(1) is immediate from the previous discussion.  For (2), the only
Gr\"obner basis elements that can be used in the reduction of $M$
are $G_1,\dots,G_{k-1}$, so the reduction process
cannot introduce a monomial divisible by any of $x_k,\dots,x_n$.
\end{proof}

One useful consequence of ``stickiness'' is the following.
\begin{cor}
\label{sticky-vs-nonsticky}
Let $\lambda$ be an indecomposable partition which is its own core,
and $k:=\lambda_1$.  Let $f=\sum_{i=1}^n a_i x_i$ be an element of 
the degree-one graded piece $R^\lambda_1$ of $R^\lambda$.
Decompose $f$ as $f=g+h$, where
$$
g~=~\sum_{i=1}^{k-1} a_i x_i, \qquad\qquad
h~=~\sum_{i=k}^n a_i x_i.
$$
%\begin{align*}
%g&=\sum_{i=1}^{k-1} a_i x_i, \\
%h&=\sum_{i=k}^n a_i x_i.
%\end{align*}
If $f^m=0$ in $R^\lambda$ for some positive integer $m$, then $g^m=0$
in $R^\lambda$.
\end{cor}
\begin{proof}
Note that $f^m = g^m + p$, where $p$ is some polynomial divisible by
$a_kx_k+\dots+a_nx_n$.  Passing to the standard forms,
we find that $0=\greduce{g^m}+\greduce{p}$.
By Corollary~\ref{sticky}, no monomial in $\greduce{g^m}$ is divisible
by a sticky variable (that is, one of $x_k,\dots,x_n$), but every
monomial in $\greduce{p}$ is divisible by a sticky variable.  Therefore
$\greduce{g^m}=0$ ($= \greduce{p}$). 
\end{proof}

%%%%%%%%%%%%%%%%%%%%%%%%%%%%%%%%%%%%%%%%%%%%%%%%%%%%%%%%%%%%%%%%%%%%%%%%%%%%%%%%%%%%%%%
\section{Nilpotence of linear forms in the cohomology of $G/B$}
\label{G/B-section}
%%%%%%%%%%%%%%%%%%%%%%%%%%%%%%%%%%%%%%%%%%%%%%%%%%%%%%%%%%%%%%%%%%%%%%%%%%%%%%%%%%%%%%%

The main result of this section, Theorem~\ref{G/B-nilpotents},
concerns the nilpotence orders of degree-$1$ elements in the graded ring 
$H^*(G/B)$.  This result may be of independent interest, and
it would be nice to have a geometric explanation for it.

Recall that $H^*(G/B) = R^{n^n} \cong \Zz[x_1,\ldots,x_n] / J$,
where
\begin{equation}
\label{G/B-cohomology-descriptions}
J =  \ideal{ e_i(n): 1 \leq i \leq n }  
  = I_{n^n} =  \ideal{ h_{n-i+1}(i): 1 \leq i \leq n }. \\
\end{equation}

We digress to discuss graded $\Zz$-algebras and nilpotence.
A {\it standard graded $\Zz$-algebra} is a ring $R$ with
a $\Zz$-module direct sum decomposition $R = \bigoplus_{d \geq 0} R_d$
in which each $R_d$ is a free $\Zz$-module, $R_d \cdot R_e \subset R_{d+e}$
and $R$ is generated over the subalgebra $R_0=\Zz$ by $R_1$.
Let $R$ be a ring and $f \in R$ a nilpotent element (that is, some power of $f$
is zero).  The \defterm{nilpotence order} of $f$
is defined as the smallest integer $k$ such that $f^k=0$; we will sometimes say that $f$ is
\defterm{$k$-nilpotent}.  (So $f$ has nilpotence order 1 if and only if $f=0$.)

By Theorem~\ref{cohomology-presentation},
$R^\lambda=\Zz[x_1,\ldots,x_n]/I_\lambda$ is a standard graded $\Zz$-algebra,
with $R^\lambda_1 \cong H^2(X_\lambda;\:\Zz)$.  Furthermore,
every element of $R^\lambda_1$ is nilpotent,
since $R^\lambda$ has finite rank as a $\Zz$-module.
The nilpotence order of these linear forms
will be our main tool in distinguishing
the rings $R^\lambda$.  In this section, we study the case that $\lambda=n^n$;
we treat the general case in Section~\ref{general-nilpotence-section}.

Note that the images of the variables $x_i$ in $R^{n^n}$ 
satisfy $x_i^n=0$.  Indeed, by Corollary~\ref{variable-symmetry},
it is sufficient to prove that $x_1^n=0$, which follows
from~\eqref{G/B-cohomology-descriptions} since $I_{n^n}$ 
contains the element $h_{n-1+1}(1)=h_n(1)=x_1^n$.  In fact, more is true:

\begin{thm} \label{G/B-nilpotents}
Let $f \in H^2(G/B) \cong (R^{n^n})_1$.

Then $f$ has nilpotence
order greater than or equal to $n$, with equality if and only if $f$ is congruent, modulo
$J$, to a scalar multiple of one of the variables $x_1,\dots,x_n$.
\end{thm}

We first show that $n$ is the minimal nilpotence order
achieved by any linear form.

\begin{prop} 
\label{smallest-G/B-nilpotence}
Let $f \in R^{n^n}_1$ be a linear form.
If $f^{n-1}=0$, then $f=0$.
\end{prop}

\begin{proof} 
Let $\hat{f}$ be a preimage of $f$ under the quotient map $\Zz[x_1,\ldots,x_n]
\twoheadrightarrow R^{n^n}$.
Then $f^{n-1}=0$ means $\hat{f}^{n-1} \in J$.
By degree considerations, this means that
$\hat{f}^{n-1}$ belongs to the ideal
    \begin{equation} \label{ideal-I}
    I := \ideal{ e_i(n): 1 \leq i \leq n-1 } ~\subset~ \Zz[x_1,\ldots,x_n].
    \end{equation}
Thus it suffices to show that $I$ is a radical ideal,
since then $\hat{f} \in I$ and $f=0$ in $R^{n^n}$.  We will show
something slightly stronger: that
the ideal
$I':=\ideal{ e_i(n): 1 \leq i \leq n-1 } \subset \Cc[x_1,\ldots,x_n]$
is radical.  Indeed,
any nonzero nilpotent in $\Zz[x_1,\ldots,x_n]/I$ would
give rise to a nonzero nilpotent in $\Cc[x_1,\ldots,x_n]/I'$.

Let $\zeta$ be a primitive $n^{th}$ root of unity.
We claim that $I'$ is the vanishing ideal $I(V)$ for the
variety $V \subset \Cc^n$, defined as the union of all lines
whose slope vector is any permutation
of $(1,\zeta,\dots,\zeta^{n-1})$.  Note that there are exactly $(n-1)!$
such lines, because two such slope vectors that differ by multiplication
by a root of unity give rise to the same line.  
Equating coefficients of powers of $t$ in the equation
$$
t^n-1=\prod_{i=1}^n (t-\zeta^i) = 
\sum_{i=0}^n e_i(1,\zeta,\dots,\zeta^{n-1}) t^i
$$
shows that $I' \subset I(V)$.  For the reverse inclusion, note
that $e_1(n),\ldots,e_n(n)$ is a regular sequence in
$\Cc[x_1,\ldots,x_n]$, and therefore cuts out scheme-theoretically
a complete intersection of Krull dimension $1$, that is,
a set of curves with various multiplicities.  By B\'ezout's Theorem,
the sum of the degrees of those curves, counted with multiplicities,
must be
$$
\deg(e_1(n)) \cdot \deg(e_2(n)) \cdots \deg(e_n(n)) ~=~ 1\cdot 2 \cdots 
(n-1) ~=~ (n-1)!
$$
But this complete intersection contains at least $(n-1)!$ lines 
in $V$, each of degree $1$.  Therefore it contains
no other curves, and each line occurs with multiplicity 1; that is, $I'=I(V)$.
\end{proof}

To complete the proof of Theorem~\ref{G/B-nilpotents},
we must show that the scalar multiples of the variables $x_i$
are the only $n$-nilpotent linear forms in $R^{n^n}$.  In what follows,
we regard a linear form $f=\sum_{i=1}^n a_i x_i$
as a $\Cc$-linear functional, mapping $v=(v_1,\dots,v_n) \in \Cc^n$
to $\sum_{i=1}^n a_i v_i$.

\begin{lem}
\label{mike's-lemma}

Let $f=\sum_{i=1}^n a_i x_i$, with $a_i \in \Cc$, and let $\alpha \in
\Cc^*$ be a nonzero constant.  Suppose that $f(v)^n=\alpha$ for all
$v \in \Cc^n$ whose coordinates are permutations of the distinct $n^{th}$
roots of unity.

Then $f \in \Cc x_i + \Cc e_1(n)$ for some $i$.
\end{lem}

\begin{proof}
Let $\zeta$ be a primitive $n^{th}$ root of unity.
Let the symmetric group $\Sym_n$ act
on $\Cc^n$ by permuting coordinates, and for a permutation $\sigma \in \Sym_n$,
abbreviate $f(\sigma(1,\zeta,\ldots, \zeta^{n-1}))$ by $f(\sigma)$.
Replacing $f$ with $f/\alpha$, we may assume 
that $f(\sigma)^n=1$ for all $\sigma \in \Sym_n$.
That $f$ has the desired form is equivalent to the statement that
at least $n-1$ of the coefficients $a_1,\ldots,a_n$ are equal.
This is trivial if $n=1$ or $n=2$, and can be checked by direct
calculation if $n=3$.
%The values
%$\{a_1+a_2\zeta+a_3\zeta^2, a_1\zeta+a_2\zeta^2+a_3, a_1\zeta^2+a_2+a_3\zeta\}$
%range over $\{1,\zeta,\zeta^2\}$.
%So do the values
%$\{a_1+a_2\zeta^2+a_3\zeta, a_1\zeta^2+a_2\zeta+a_3, a_1\zeta+a_2+a_3\zeta^2\}$.
%Find the guy from each set that equals 1.  All of the nine possibilities work
%the same way; suppose that $a_1+a_2\zeta+a_3\zeta^2=a_1+a_2\zeta^2+a_3\zeta=1$.
%Then $a_2(\zeta-\zeta^2)-a_3(\zeta^2-\zeta)=0$, i.e., $a_2=a_3$ as desired.
Therefore, suppose $n \geq 4$.  By transitivity, it suffices to show that
if two coefficients $a_i$ are different, then the other $n-2$ are mutually equal.

Suppose that $a_1 \neq a_2$. Choose $i\neq j\in [n]$ so as to
maximize $|\zeta^i-\zeta^j|$, and let $\sigma \in \Sym_n$ such that
$\sigma(1) = i$ and $\sigma(2) = j$.  Then $f(\sigma)$ and
$f((12)\circ\sigma)$ are both $n^{th}$ roots of unity, and
  \begin{equation} \label{nth-roots-eqn}
    f(\sigma)-f((12)\circ\sigma) = (a_1-a_2)(\zeta^i-\zeta^j).
  \end{equation}
Taking the magnitude of both sides, the
choice of $i$ and $j$ implies that $|a_1-a_2| \leq 1$.
On the other hand, if we choose $i'\neq j'\in [n]$ to
minimize $|\zeta^{i'}-\zeta^{j'}|$, the same argument implies
that $|a_1-a_2| \geq 1$.  We conclude that $|a_1-a_2| = 1$.

Note that $\zeta^i$ and $\zeta^j$ are the only $n^{th}$ roots of unity whose
difference is $\zeta^i-\zeta^j$.
(This may be seen most easily by plotting the $n^{th}$ roots of unity in
the complex plane, and observing that no two of the line segments
joining two maximally distant roots are parallel.)
Therefore, the equation~\eqref{nth-roots-eqn} implies that
the values $f(\sigma)$ and $f((12)\circ\sigma)$ do not depend on
$\sigma(3),\dots,\sigma(n)$.  Hence $a_3 = \ldots = a_n$ as desired.
\end{proof}

\begin{prop} \label{G/B-minimal-nilpotents}
Let $f \in R^{n^n}_1$ be a linear form such that $f^n=0$.

Then $f \in \Zz x_i$ for some $i$.
\end{prop}

\begin{proof}
Let $\hat{f}$ be a preimage of $f$ under the quotient map $\Zz[x_1,\ldots,x_n]
\twoheadrightarrow R^{n^n}$;
that is, $\hat{f}^n \in J$.  By degree considerations,
there is a constant $\alpha \in \Zz$ such that
$\hat{f}^n \equiv \alpha e_n(n)$ modulo $I$.
As in the proof of Proposition~\ref{smallest-G/B-nilpotence}, the ideal
$I$ vanishes on all vectors $v$ whose
coordinates are a permutation of the distinct $n^{th}$ roots of unity.  Therefore
$\hat{f}^n(v) = \alpha e_n(n)(v) = (-1)^{n-1} \alpha$
for all such vectors $v$.  By Lemma~\ref{mike's-lemma},
there is some $i$ such that $\hat{f} \in \Cc x_i + \Cc e_1(n)$.  As 
$\hat{f} \in \Zz[x_1,\ldots,x_n]$, this implies 
$\hat{f} \in \Zz x_i + \Zz e_1(n)$.  Consequently
$f \in \Zz x_i$ in $R^{n^n}$.  This completes the proof of the proposition
and of Theorem~\ref{G/B-nilpotents}.
\end{proof}

%%%%%%%%%%%%%%%%%%%%%%%%%%%%%%%%%%%%%%%%%%%%%%%%%%%%%%%%%%%%%%%%%%%%%%%%%%%%%%%%%%%%%%%
\section{Nilpotence of linear forms in the cohomology of $X_\lambda$}
\label{general-nilpotence-section}
%%%%%%%%%%%%%%%%%%%%%%%%%%%%%%%%%%%%%%%%%%%%%%%%%%%%%%%%%%%%%%%%%%%%%%%%%%%%%%%%%%%%%%%

Throughout this section, $\lambda$ will be an indecomposable partition.
We continue our study of nilpotence orders of
linear forms in the graded $\Zz$-algebra
$R^\lambda=H^*(X^\lambda)$.  The main result is the following
classification of linear forms of minimal nilpotence order, generalizing
Theorem~\ref{G/B-nilpotents}. 

\begin{thm} \label{classification-of-nilpotents}
Let $\lambda = (0 < \lambda_1 \leq \cdots \leq \lambda_n)$ 
be an indecomposable partition, and let $k:=\lambda_1$.
Then $k$ is the minimal nilpotence order of any linear form
in $R^\lambda$.  Moreover, if $\lambda$ has exactly $m$ parts
equal to $k$ (that is, $k = \lambda_1 = \cdots = \lambda_m < \lambda_{m+1}$),
then the elements of $R^\lambda_1$ of nilpotence order exactly
$k$ are classified as follows:

\begin{enumerate}
\item[\sf Case I.] Either $\lambda_{k-1}>k$, or $n < k-1$.

  Then the $k$-nilpotents in $R^\lambda_1$ are the multiples of $x_1,\ldots,x_m$.

\item[]

\item[\sf Case II.] $\lambda_{k-1}=k$ (that is, $m=k-1$).

 \begin{enumerate}
 \item[\sf Subcase IIa.] Either $\lambda_k > k+1$, or $k$ is odd.

  Then the $k$-nilpotents are $x_1,\dots,x_{k-1}$, and $ x_1+\cdots+x_{k-1}$.
 \item[]
 \item[\sf Subcase IIb.] Both $\lambda_k=k+1$ and  $k$ is even.

  Then the $k$-nilpotents are 
  $x_1,\dots, x_{k-1}$, $x_1+\cdots+x_{k-1}$, and $x_1+\cdots+x_{k-1}+2x_k$.
 \end{enumerate}
\end{enumerate}
\end{thm}

By way of motivation for the rather technical matter of this section,
we explain how the classification of nilpotents will be used in the next
two sections to
recover a partition from its cohomology ring.  Theorem~\ref{classification-of-nilpotents}
implies immediately that $\lambda_1$ is an isomorphism invariant of $R^\lambda$.
Moreover, by the presentation of Theorem~\ref{cohomology-presentation},
the quotient ring $R^\lambda\:/\ideal{x_1}$
may be identified with the ring $R^\nu$, where
$\nu = (\lambda_2-1, \lambda_3-1\ldots, \lambda_n-1)$
is the partition obtained by removing the first row and column from $\lambda$
(see Figure~\ref{nu-example}).
However, it is really necessary to describe $R^\nu$ as a quotient
$R^\lambda\:/\ideal{f}$, where $f$ is some linear form identified
intrinsically from the structure of $R^\lambda$ as a standard graded $\Zz$-algebra,
that is, in a way that 
\emph{does not depend on the presentation}.  The classification
of nilpotents in Theorem~\ref{classification-of-nilpotents} is the tool 
that allows this.   It turns out
that we will require almost all, but not quite all of the last assertion
in the theorem, so we only prove the parts that will be used.
(The arguments we omit are very similar to those that we include.)

\vskip .2in
In the first part of this section, culminating in Proposition~\ref{min-nil-deg},
we prove the first assertion of Theorem~\ref{classification-of-nilpotents}, namely
that $k = \lambda_1$ is
the minimal nilpotence order of any linear form in $R^\lambda$.  We begin with
a weaker statement: that no linear form in the first $k-1$ variables
has nilpotence order less than $k$.

\begin{lem} \label{no-initial-nilpotent}
Let $\lambda$ be indecomposable with $k:=\lambda_1$.
Let $f = \sum_{i=1}^{k-1} a_i x_i \in R^\lambda_1$; that is,
$f$ is supported only on the first $k-1$ variables.  Then, in $R^\lambda$,

\begin{enumerate}
\item[(a)] $f^{k-1}=0$ if and only if $f=0$, and
\item[(b)] if $f^k=0$, then $f$ is a scalar multiple of one of the following:
$x_1, \ \ldots, \ x_{k-1}, \ x_1+\cdots+x_{k-1}$.
\end{enumerate}
\end{lem}

\begin{proof}
By Proposition~\ref{Leray-Hirsch} and the hypothesis that 
$f$ is supported only on the first $k-1$ variables, we may assume without
loss of generality that $n \leq k-1$.
By Corollary~\ref{contain-epi}, we may decrease the part sizes of $\lambda$ 
(if necessary), so as to assume that $\lambda = k^n$.
But then using Proposition~\ref{Leray-Hirsch} again, we 
can re-introduce parts $\lambda_{n+1},\lambda_{n+2},\ldots,\lambda_k$
all of size $k$, and work in the ring $R^{k^k} \cong GL_k(\Cc)/B$, where
assertion (a) follows from Theorem~\ref{G/B-nilpotents}.  

In fact, assertion (b) also follows from 
Theorem~\ref{G/B-nilpotents}.  The degree-$1$ graded piece of
$I_{k^k}$ is generated by $e_1(k)$,
so the elements of $R^{k^k}_1$ listed above are the only ones that are congruent
modulo $I_{k^k}$ to a scalar multiple of a variable $x_i$ (here
we use the fact that $x_1 + \cdots + x_{k-1} = e_1(k) - x_k$).
\end{proof}

An immediate consequence of Lemma~\ref{no-initial-nilpotent}
is that every linear form of nilpotence order $\lambda_1-1$ must
be supported on at least one of the variables $x_k, \ldots, x_n$.
This is where the concept of ``stickiness'' introduced in
Corollary~\ref{sticky} first comes into play.
% In what follows, we work with the free $\Zz$-basis for $R^\lambda$
% given by the standard monomials.  Thus, if $f \in \Zz[x_1,\dots,x_n]$,
%We will use the following notation: if $f$ is a polynomial
%and $M$ is a standard monomial, then $[M]f$ denotes the coefficient of
%$M$ in $f$.

\begin{prop}\label{no-k-1}
Let $\lambda$ be indecomposable with $k:=\lambda_1$, and
let $f \in R^\lambda_1$.
Then $f^{k-1}=0$ if and only if $f=0$ in $R^\lambda$.
\end{prop}

\begin{proof}
Assume $f \neq 0 \in R^\lambda_1$, but $f^{k-1}=0$ in $R^\lambda$.
By Lemma~\ref{no-initial-nilpotent}(a), we may assume $n \geq k$.
By Proposition~\ref{Leray-Hirsch}, we may assume without
loss of generality that $\lambda$ is its own core.

Writing $f=g+h$, where $g = a_1x_1+\dots+a_{k-1}x_{k-1}$
and $h=a_k x_k +\dots+a_n x_n$, it follows from 
Corollary~\ref{sticky-vs-nonsticky} that
$g^{k-1}=0$.  Hence $g=0$ by 
Lemma~\ref{no-initial-nilpotent}.  That is, $f=h$.
If $f$ is not supported on $x_n$ (that is, $a_n=0$), then we may replace
$\lambda$ with the partition obtained by removing the $n$th (largest) row.
Repeating this as many times as necessary, we may assume without
loss of generality that $a_n \neq 0$.

Now let $M$ be any monomial in the variables $x_1,\ldots,x_{k-1}$.
Note that
    \begin{equation} \label{sticky1}
    \coef{x_n M}{\greduce{f^{k-1}}} ~=~ \coef{x_n M}{\greduce{(a_nx_n)^{k-1}}}
    \end{equation}
because the variables $x_k,\dots,x_{n-1}$ are sticky (Corollary~\ref{sticky}).  
Reducing $(a_nx_n)^{k-1}$ using the Gr\"obner basis element $G_n$
of~\eqref{RGB}, we find that
%$x_n^2 \equiv -x_n(x_1+\dots+x_{n-1})$ modulo $I_\lambda$.  Thus, in $R^\lambda$, we have
    \begin{equation} \label{repeated-reduction}
    \begin{aligned}
    (a_nx_n)^{k-1} ~&=~ -a_nx_n^{k-2}(x_1+\dots+x_{n-1})\\
    &=~ a_n^2x_n^{k-3}(x_1+\dots+x_{n-1})^2\\
    &=~ \dots\\
    &=~ \alpha x_n (x_1+\dots+x_{n-1})^{k-2},
    \end{aligned}
    \end{equation}
%
%
%, and
%can use this to do reduction $k-2$ times and conclude
%    $$(a_nx_n)^{k-1} ~=~ \alpha x_n (x_1+\dots+x_{n-1})^{k-2}$$
%in $R^\lambda$, 
where $\alpha=(-1)^{k-2} a_n^{k-2} \neq 0$.
Combining this with~\eqref{sticky1} yields
    \begin{subequations}
    \begin{eqnarray}
    \coef{x_n M}{\greduce{f^{k-1}}}
    &=& \alpha \coef{x_n M}{\greduce{x_n (x_1+\dots+x_{n-1})^{k-2}}} \notag\\
    &=& \alpha \coef{x_n M}{\greduce{x_n (x_1+\dots+x_{k-1})^{k-2}}} \label{dirty-trick}\\
    &=& \alpha \coef{M}{\greduce{(x_1+\dots+x_{k-1})^{k-2}}} \label{final-dirty-trick}
    \end{eqnarray}
    \end{subequations}
where~\eqref{dirty-trick} follows from stickiness, and~\eqref{final-dirty-trick}
from the fact that only $G_1,\dots,G_{k-1}$ are used in reducing~\eqref{dirty-trick}.

The polynomial $x_1+\dots+x_{k-1}$ is nonzero in $R^\lambda$ since $\lambda$
is indecomposable.  Thus Lemma~\ref{no-initial-nilpotent} implies that
$(x_1+\dots+x_{k-1})^{k-2} \neq 0$ as well, and so 
there exists some monomial $M$ in the variables $x_1,\ldots,x_{k-1}$
for which $[M]\greduce{(x_1+\dots+x_{k-1})^{k-2}} \neq 0$.  
Note that $x_n M$ is also a standard monomial for $I_\lambda$.
Therefore $[x_n M]\greduce{f^{k-1}} \neq 0$, a contradiction.
\end{proof}

\begin{prop} \label{min-nil-deg}
When $\lambda$ is indecomposable, the 
number $k=\lambda_1$ is an isomorphism invariant of $R^\lambda$ as a graded ring:
namely, it is the minimum nilpotence order achieved by any linear form.
\end{prop}

\begin{proof}
Proposition~\ref{no-k-1} states that no nonzero linear form can have nilpotence
order strictly less than $k=\lambda_1$.  On the other hand, $x_1$ has nilpotence
order at most $k$, because $x_1^k = h_k(1) \in I_\lambda$.
\end{proof}

In the second part of this section, we show that the various 
linear forms mentioned in Theorem~\ref{classification-of-nilpotents}
are the only possible $k$-nilpotents in $R^\lambda$.
We begin by determining the nilpotence order of each variable.

\begin{prop}
\label{variable-nilpotence-indecomp}
When $\lambda$ is indecomposable, the variable $x_i$ is $\lambda_i$-nilpotent
in $R^\lambda$.
\end{prop}

\begin{proof}
Let $p = \lambda_i$.  First, we show that $x_i^p = 0$ in $R^\lambda$.  Let
$\kappa$ be the partition given by
    $$\kappa ~:=~ (\underbrace{p,\dots,p}_{i\text{ times}},\;
    \lambda_{i+1},\; \lambda_{i+2},\; \ldots,\; \lambda_n).$$
Then $\lambda$ is a subpartition of $\kappa$,
so $R^\lambda$ is a quotient of $R^\kappa$ by Lemma~\ref{contain-epi}.
It suffices to show that $x_i^p = 0$ in $R^\kappa$, which
follows from Corollary~\ref{variable-symmetry} since
$x_1^p \in I_\kappa$.

It remains to show that $x_i^{p-1} \neq 0$ in $R^\lambda$.
By Proposition~\ref{Leray-Hirsch} and Corollary~\ref{contain-epi},
it suffices to show that $x_i^{p-1} \neq 0$ in $R^\mu$,
where $\mu$ is the subpartition of $\lambda$ given by
    $$\mu ~:=~ (2,3,\dots,i-1,i,p).$$
Note that $\mu$ is indecomposable, and that $R^\mu$ has a reduced Gr\"obner
basis given by~\eqref{putative-GB}.
A Gr\"obner reduction similar to~\eqref{repeated-reduction},
using the Gr\"obner basis element $x_i^{p-i+1}+x_i^{p-i}h_1(i-1)$
yields the equation
    $$x_i^{p-1} ~\equiv~ (-1)^{i-1} x_i^{p-i} h_1(i-1)^{i-1} \pmod{I_\mu}.$$
Since further reductions modulo $I_\mu$
can only involve the other generators $h_2(1),h_2(2),\ldots,h_2(i-1)$,
we may conclude that $x_i^{p-1} \neq 0$ in $R^\mu$, provided that $h_1(i-1)^{i-1} \neq 0$
in $R^{(2,3,\cdots,i-1,i)}$.  
Using the fact that $h_1(i-1)=e_1(i-1)$, this follows from the 
following more general assertion: for any $m \geq 1$ and $i \geq 1$, 
    \begin{equation}
    \label{power-to-index-assertion}
    e_1(i-1)^m ~=~ e_m(i-1) ~\neq~ 0 \qquad \text{in } R^{(2,3,\cdots,i-1,i)}.
    \end{equation}
This is trivially true for $i\leq 2$.  For $i>2$, we prove it by induction on $i$:
%Assertion \eqref{power-to-index-assertion}
%follows by induction on $i$:
  \begin{align*}
  e_1(i-1)^m ~&~= (x_{i-1} + e_1(i-2))^m \\
              &~= \sum_{j=0}^m  \binom{m}{j} x_{i-1}^j e_1(i-2)^{m-j}\\
              &~= e_1(i-2)^m + \sum_{j=1}^m  \binom{m}{j} x_{i-1}^j e_1(i-2)^{m-j}\\
              &~\equiv e_1(i-2)^m +
                 \sum_{j=1}^m \binom{m}{j} (-1)^{j-1} x_{i-1} e_1(i-2)^{m-1}
                 \pmod{I^{(2,3,\cdots,i-1,i)}}.
  \end{align*}
This last expression follows from using $x_{i-1}h_1(i-1) = x_{i-1}^2+x_1h_1(i-2)
= x_{i-1}^2+x_1h_1(i-2)$ to perform repeated Gr\"obner reduction on each summand.
By induction, $e_1(i-2)^m = e_m(i-2)$, so we obtain
  \begin{align*}
  e_1(i-1)^m ~&=~ e_m(i-2) +
                   x_{i-1} e_{m-1}(i-2) \sum_{j=1}^m \binom{m}{j} (-1)^{j-1} \\
              &=~ e_m(i-2) + x_{i-1} e_{m-1}(i-2)  \\
              &=~ e_m(i-1),
  \end{align*}
establishing~\eqref{power-to-index-assertion} as desired.
\end{proof}

\begin{prop}
Let $f=\sum_{i=1}^n a_i x_i \in R^\lambda$.  Suppose that $f^k=0$.

Then $f$ is a scalar multiple of one of the following:
    \begin{equation} \label{candidate-nilpotents}
    \begin{aligned}
    & x_1,\, x_2,\, \dots,\, x_{k-1}, \\
    & x_1+\cdots+x_{k-1}, \\
    & x_1+\cdots+x_{k-1}+2x_k.
    \end{aligned}
    \end{equation}
The last case can occur only if $k$ is even.
\end{prop}

\begin{proof}
By Corollary~\ref{contain-epi}, we may replace $\lambda$ with its core. 
Let $g=\sum_{i=1}^{k-1} a_i x_i$ be the part of $f$ in the non-sticky variables.
Then $g^k=0$ by Corollary~\ref{sticky-vs-nonsticky}.
By Lemma~\ref{no-initial-nilpotent}(b), $g$ is either zero or of the form
$\alpha x_i$ for some $i \in \{1,2,\dots,k-1\}$, or
$\alpha(x_1+\dots+x_{k-1})$, where $\alpha$ is a nonzero scalar.
Without loss of generality, we may assume that $\alpha=1$.

If $f=g$ then we are done.
Otherwise, we must show that $f$ is a scalar multiple of 
$x_1+\cdots+x_{k-1}+2x_k$, and $k$ is even.
By Proposition~\ref{Leray-Hirsch},
we may assume without loss of generality that $f$ involves the
variable $x_n$ with non-zero coefficient; that is, 
    $$f = g + h + ax_n,$$
where $a:=a_n \neq 0$ and 
$h$ is a linear form in the variables $x_k,\dots,x_{n-1}$.
We consider in turn each of the three possibilities: namely, $g=0$, $g=x_i$,
or $g=x_1+\dots+x_{k-1}$.

\vskip 0.1in
\noindent{\sf Case 1:} $g=0$.

We will rule out this case by deriving a contradiction
from the assumption that $f^k=0$ in $R^\lambda$.
Taking the further quotient of $R^\lambda$ by the variables $x_k,\dots,x_{n-1}$,
one obtains a ring isomorphic to $R^\mu$, where 
$$
\mu=(\underbrace{k,\ldots,k}_{k-1\text{ times}},k+1)
$$
is an indecomposable partition, with $k$ parts, equal to its own core.
If $f^k=0$ in $R^\lambda$,
then $(ax_k)^k = a^k x_k^k = 0$ in $R^\mu$.  So $x_k^k=0$ in $R^\mu$ (because
$a \neq 0$).  But this contradicts Corollary~\ref{variable-nilpotence-indecomp},
since $\mu_k=k+1$.

\vskip 0.1in
\noindent{\sf Case 2:} $g= x_i$, where $i \in \{1,2,\dots,k-1\}$.

Assume that $k\geq 3$ (the case $k=2$ falls under Case~3 below).
As in Case 1, we wish to reach a contradiction.
Consider the quotient ring
    $$S ~:=~ R^\lambda / \ideal{x_k,~\dots,~x_{n-1},~x_1+x_2+\dots+x_{k-1}+x_n},$$
which is isomorphic to $R^{k^k}$.  Let $\tilde f = x_i - a(x_1+\dots+x_{k-1})$
be the image of $f$ in $S$; then $\tilde f^k=0$.
By Theorem~\ref{G/B-nilpotents},
$\tilde f$ must be a scalar multiple of some variable.  This is possible only if $k=3$
and $a=1$; that is, $f$ is a scalar multiple of either $x_1+x_3$
or $x_2+x_3$.  All that remains is to check that neither
$(x_1+x_3)^3$ nor $(x_2+x_3)^3$ belongs to the ideal
$I_{3^3}=\ideal{ h_3(1), h_2(2), h_2(3) }$; this is a routine calculation.
Thus $f^k \neq 0$ in all cases, a contradiction.
Case~2 is therefore ruled out.

\vskip 0.1in
\noindent{\sf Case 3:} $g = x_1+\ldots+x_{k-1}$.

Let $M$ be any standard monomial for $I_\lambda$
of degree $k-1$ in the non-sticky variables
$x_1,\dots,x_{k-1}$; then $x_n M$ is also standard.  Using stickiness
of the variables $x_k,\ldots,x_{n-1}$ 
and the fact that $G_n = x_n(x_1+\dots+x_n) \in I_\lambda$,
we have for every such monomial
 \begin{align*}
    [x_n M] f^k ~&=~ [x_n M] (g+ax_n)^k 
    ~=~ [x_n M] \sum_{i=0}^k \binom{k}{i} a^i x_n^i g^{k-i} 
    ~=~ [x_n M] \sum_{i=1}^k \binom{k}{i} a^i x_n^i g^{k-i} \\
    &=~ [x_n M] \sum_{i=1}^k 
         \binom{k}{i} a^i g^{k-i} (-1)^{i-1} x_n(x_1+\dots+x_{n-1})^{i-1} \\
    &=~ [x_n M] \sum_{i=1}^k 
         \binom{k}{i} a^i g^{k-i} (-1)^{i-1} x_n(x_1+\dots+x_{k-1})^{i-1} \\
    &=~ [M] \sum_{i=1}^k \binom{k}{i} a^i g^{k-1} (-1)^{i-1} 
    ~=~ \left(\sum_{i=1}^k \binom{k}{i} a^i (-1)^{i-1} \right) [M]g^{k-1} \\
    &=~ \left(1-(1-a)^k\right) [M]g^{k-1}.
 \end{align*}
This last expression must be zero since $f^k=0$ in $R^\lambda$.  On the other hand,
$g^{k-1}\neq 0$ in $R^\lambda$, 
so there is at least one such monomial $M$ in $x_1,\ldots,x_{k-1}$
for which $[M]g^{k-1}\neq 0$.
It follows that $1-(1-a)^k=0$.  Since $a \neq 0$, the only possibility is that $k$ is even
and $a=2$.  If $n=k$, then we are done; we need to rule out the case $n>k$.

Suppose that $n>k$.  Replacing $x_n$ with $x_k$ in the above calculation,
we find that the coefficient $a_k$ is either 0 or 2.  Bearing in mind
that $g+x_k = x_1+\cdots+x_{k-1}+x_k=h_1(k)$, we pass to the quotient ring
  \begin{eqnarray*}
    T &:=& R^\lambda\:/\ideal{x_{k+1},\, x_{k+2},\, \dots,\, x_{n-1},\, g+x_k}\\
      &\cong& \Zz[x_1,\dots,x_k,x_n]\:/
           \ideal{h_k(1),\, h_{k-1}(2),\, \dots,\, h_2(k-1) ,\,x_k(g+x_k),\,
          x_n(g+x_k+x_n),\, g+x_k}\\
      &\cong& \Zz[x_1,\dots,x_k,x_n]\:/
            \ideal{h_k(1),\, h_{k-1}(2),\, \dots,\, h_2(k-1),\,
          g+x_k,\, x_n^2}\\
      &\cong& R^{k^k} \,[x_n]\:/\ideal{x_n^2}.
  \end{eqnarray*}
Note that since $f$ equals either $g+2x_n$ or $g+2x_k+2x_n$,
and $x_k = -g$ in $T$, the image $p$ of $f$ in $T$ is of the form
$p=\pm g + 2x_n$.
Since $x_n^2$ and $g^k$ are both zero in $T$, we have
    $$p^k ~=~ \sum_{j=0}^k \binom{k}{j} (2x_n)^j (\pm g)^{k-j} ~=~ \pm 2k x_n g^{k-1}.$$
But $g^{k-1} \neq 0$ in $R^{k^k}$ by Theorem~\ref{G/B-nilpotents},
so $x_n g^{k-1} \neq 0$ in $T$.  Hence 
$p^k \neq 0$ in $T$, which implies that $f^k \neq 0$ in $R^\lambda$, as desired.
\end{proof}

We now know that every $k$-nilpotent linear form in $R^\lambda$ is,
up to scalar multiplication, one of the linear forms~\eqref{candidate-nilpotents}.
However, if $\lambda$ is not its own core,
then we must consider the possibility that one or more of these linear forms
actually has nilpotence order
strictly greater than $k$.  We examine each candidate in turn;
Proposition~\ref{variable-nilpotence-indecomp} immediately takes care of 
the possible nilpotents $x_1,\dots,x_{k-1}$.

\begin{prop} \label{nosum}
Let $\lambda$ be indecomposable with $n \geq k-1$ parts and $k=\lambda_1$.
Let $g=x_1+\cdots+x_{k-1} \in R^\lambda$.

Then $g^k=0$ if and only if $\lambda_1 = \cdots = \lambda_{k-1} = k$. 
\end{prop}

\begin{proof}
By Proposition~\ref{Leray-Hirsch}, we may assume
that $n=k-1$.  Suppose that $\lambda_1 = \cdots = \lambda_{k-1} = k$.  Then
  \begin{eqnarray*}
  R^\lambda = R^{k^{k-1}} 
    &=& \Zz[x_1,\dots,x_{k-1}]\phantom{,x_k} ~/~
          \ideal{h_k(1),\, h_{k-1}(2),\, \ldots,\,h_2(k-1)} \\
    &\cong& \Zz[x_1,\dots,x_{k-1},x_k] ~/~
          \ideal{h_k(1),\, h_{k-1}(2),\, \ldots,\,h_2(k-1),h_1(k)}\\
    &=& R^{k^k},\\
  \end{eqnarray*}
and $g = -x_k$ in $R^{k^k}$, so $g^k=0$ follows
from Theorem~\ref{G/B-nilpotents}.

Conversely, suppose that $\lambda_{k-1}>k$.  We will show that $g^k\neq 0$.
Let $\mu$ be the subpartition of $\lambda$ given by
$$
\mu=(\underbrace{k-1,\ldots,k-1}_{k-2 \text{ times}},k+1)
$$ 
(see Figure~\ref{lambda-example4}).
\begin{figure}
\begin{center}
\resizebox{4.8cm}{2.4cm}{\includegraphics{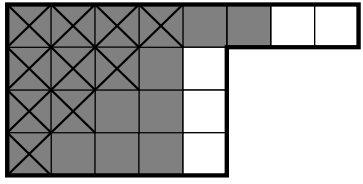}}  % revised JLM 2/26
\end{center}
\caption{\label{lambda-example4}
The subpartition $\mu$ of Proposition~\ref{nosum} (shaded).}
\end{figure}
By Corollary~\ref{contain-epi},
it will suffice to show that $g^k \neq 0$ in $R^\mu$.  We may rewrite the
presentation of $R^\mu$ as
  \begin{eqnarray*}
  R^\mu
    &=& \Zz[x_1,\dots,x_{k-1}]\:/\ideal{h_{k-1}(1),\, h_{k-2}(2),\,
      \dots,\, h_2(k-2),\, h_3(k-1)} \\
    &=& \Zz[x_1,\dots,x_{k-1}]\:/\ideal{h_{k-1}(1),\, h_{k-2}(2),\,
      \dots,\, h_2(k-2),\, x_{k-1}^3 + x_{k-1}^2h_1(k-2)},
  \end{eqnarray*}
using the fact that 
  \begin{eqnarray*}
    h_3(k-1) &=& x_{k-1}^3 + x_{k-1}^2h_1(k-2) + x_{k-1}h_2(k-2) + h_3(k-2)\\
    &=& x_{k-1}^3 + x_{k-1}^2h_1(k-2) + x_{k-1}h_2(k-2) + h_3(k-3) + x_{k-2}h_2(k-2),\\
  \end{eqnarray*}
Therefore $x_{k-1}^j \equiv (-1)^j x_{k-1}^2 h^{j-2}$ for all $j \geq 3$.
Letting $h=h_1(k-2)=x_1+\dots+x_{k-2}$, so that $g=h+x_{k-1}$, we have in $R^\mu$
  \begin{eqnarray*}
    g^k 
% &=& \sum_{j=0}^k \binom{k}{j} h^{k-j} x_{k-1}^j \\
    &=& h^k + \binom{k}{1} h^{k-1} x_{k-1} + \binom{k}{2} h^{k-2} x_{k-1}^2 +
      \sum_{j=3}^k \binom{k}{j} h^{k-j} (-1)^j x_{k-1}^2 h^{j-2} \\
    &=& h^k + k h^{k-1} x_{k-1} + h^{k-2} x_{k-1}^2\sum_{j=2}^k (-1)^j \binom{k}{j} \\
    &=& h^k + k h^{k-1} x_{k-1} + (k-1) h^{k-2} x_{k-1}^2.
  \end{eqnarray*}
No further Gr\"obner reduction is possible, so
$g^k$ is zero if and only if $h^k$, $k h^{k-1}$, and $(k-1) h^{k-2}$ are all zero.
But $k>1$, and $h^{k-2} \neq 0$ by Proposition~\ref{no-k-1}.
We conclude that $g^k \neq 0$ in $R^\mu$ as desired.
\end{proof}

For the remaining assertions of Theorem~\ref{classification-of-nilpotents},
we are left only to consider the potentially $k$-nilpotent linear form
$g=x_1 + \ldots + x_{k-1} + 2x_k$.
Rather than determining exactly when $g$ is $k$-nilpotent as in the theorem
(which can be done by an argument similar to Proposition~\ref{nosum}),
we content ourselves with checking directly the case $k=2$, since
this is all we need for the present study.  Here $g=x_1+2x_2$,
and by Proposition~\ref{Leray-Hirsch}) we may work in the ring
$$
R^{(2,\lambda_2)} = \Zz[x_1,x_2] / \ideal{x_1^2, h_{\lambda_2-1}(2)}.
$$
Then it is easily seen that $g^2=x_1^2+4x_1 x_2 + 4x^2$ is zero in this ring
if and only if $\lambda_2=3$.

%\begin{prop} Let $\lambda$ be indecomposable, with $k=\lambda_1$ even.
%Let $g=x_1+\dots+x_{k-1}+2x_k \in R^\lambda$.
%
%Then $g$ is nilpotent of order $k$ if and only if 
%$\lambda_1=\cdots=\lambda_{k-1} = k$ and $\lambda_k = k+1$.
%\end{prop}
%
%\begin{proof}
%Suppose that $\lambda$ satisfies these conditions.  In particular,
%$\lambda$ is its own core, so~\eqref{RGB} is a reduced Gr\"obner basis
%for $I_\lambda$.  We can then check that $g^k \in I_\lambda$ by a
%calculation similar to those in the previous two propositions.
%
%Letting $h=x_1+\dots+x_{k-1}$, we have in $R^\lambda$
%  \begin{eqnarray*}
%    g^k &=& \sum_{j=0}^k \binom{k}{j} 2^j x_k^j h^{k-j}\\
%    &=& h^k + 2k x_k h^{k-1} + 
%\sum_{j=2}^k \binom{k}{j} 2^j ((-1)^{j-1} x_k h^{j-1}) h^{k-j}\\
%    &=& h^k + 2k x_k h^{k-1} - \sum_{j=2}^k \binom{k}{j} (-2)^j x_k h^{k-1}\\
%    &=& h^k - x_k h^{k-1} \sum_{j=2}^k \binom{k}{j} (-2)^j \\
%    &=& h^k - x_k h^{k-1} ((-1)^k-1)\\
%    &=& h^k = 0.
%  \end{eqnarray*}
%
%
%Second, we show the necessity.
%The variable $x_k$ is sticky, so if $g$ is nilpotent then so is $x_1+\dots+x_{k-1}$.
%Therefore, the conditions $\lambda_1=\cdots=\lambda_{k-1} = k$ is necessary by
%Proposition~\ref{no-k-1}.  If these conditions hold but $\lambda_k>k+1$,
%then a calculation similar to that of Proposition~\ref{nosum} shows that
%$g^k \neq 0$.
%\end{proof}

%%%%%%%%%%%%%%%%%%%%%%%%%%%%%%%%%%%%%%%%%%%%%%%%%%%%%%%%%%%%%%%%%%%%%%%%%%%%%%%%%%%%%%%
\section{The indecomposable case}
\label{indecomposable}
%%%%%%%%%%%%%%%%%%%%%%%%%%%%%%%%%%%%%%%%%%%%%%%%%%%%%%%%%%%%%%%%%%%%%%%%%%%%%%%%%%%%%%%

We now use the results of the previous section to prove that
an indecomposable partition is determined uniquely by the cohomology ring
of the corresponding Schubert variety.

\begin{thm}\label{mainthm-indecomp}
Every indecomposable partition $\lambda$ may be recovered from the structure
of the ring $R^\lambda$ as a graded $\Zz$-algebra.  In particular, if
$\lambda$ and $\mu$ are different indecomposable partitions, then $R^\lambda$
and $R^\mu$ are not isomorphic.
\end{thm}

\begin{proof}
We induct on $n$, the number of parts of $\lambda$.
Since $\lambda$ is indecomposable, $n$
is the rank of $R^\lambda_1$ as a free $\Zz$-module.
By Theorem~\ref{classification-of-nilpotents},
the smallest part $k:=\lambda_1$ is the minimal nilpotence order of any
member of $R^\lambda_1$.  Moreover, as mentioned at the
beginning of Section~\ref{general-nilpotence-section}, $R^\lambda/\ideal{x_1}
\cong R^\nu$, where $\nu$
is obtained from $\lambda$ by deleting the first row and
column (see Figure~\ref{nu-example}).
By induction, it suffices to show that we can describe $R^\nu$ up to isomorphism
in a way that is independent of the presentation.  

We proceed by examining the same two cases as in
Theorem~\ref{classification-of-nilpotents}; however, we subdivide Case II
slightly differently into subcases.

\vskip 0.1in
\noindent
{\sf Case I.} $\lambda_{k-1}>k$ or $n < k-1$.

Let $m$ be the greatest index such that $\lambda_m=k$.
Then Theorem~\ref{classification-of-nilpotents} tells us that
the $k$-nilpotent linear forms in $(R^\lambda)_1$ are 
(up to $\Zz$-multiples) $x_1,\ldots, x_m$.  Consequently, up to
sign, these are exactly the {\it  primitive}
$k$-nilpotents, that is, those $k$-nilpotents $f$ which
can only be expressed as a scalar multiple $\alpha g$ for another $k$-nilpotent $g$
and $\alpha \in \Zz$ if $\alpha = \pm 1$.

By Corollary~\ref{variable-symmetry}, 
one has $R^\lambda/\ideal{x_i} \cong R^\lambda/\ideal{x_1}$ ($\cong R^\nu$)
for every $i \in \{1,2,\dots,m\}$, and
hence $R^\nu$ may be identified intrinsically as the quotient of 
$R^\lambda$ by an arbitrary primitive $k$-nilpotent linear form.

\vskip 0.1in
\noindent
{\sf Case II.} $\lambda_{k-1}=k$.

Then the primitive $k$-nilpotents are (up to sign)
$x_1,\; \ldots,\; x_{k-1},\; x_1+\cdots+x_{k_1}$, and
if $k$ is even, possibly also $x_1+\cdots+x_{k-1}+2x_k$.

\vskip 0.1in
\noindent
{\sf Subcase IIA.} $k > 2$.

If $k$ is odd, then the ``extraneous''
primitive $k$-nilpotent  $x_1+\cdots+x_{k-1}+2x_k$ is absent.  If 
$k$ is even, then $x_1+\cdots+x_{k-1}+2x_k$ is distinguished
intrinsically as the unique primitive $k$-nilpotent
which is $\Zz$-linearly independent of all the others.
%we can distinguish 
%$x_1+\cdots+x_{k-1}+2x_k$ intrinsically from the rest of the 
%primitive $k$-nilpotents because its $\Zz$-multiples are the only $k$-nilpotents
%which are not expressible as $\Zz$-combinations of primitive $k$-nilpotents
%from the $k$ other $\pm$-equivalence classes.

Thus, in all cases when $k > 2$, we can intrinsically identify the
primitive $k$-nilpotents $x_1$, $\ldots$, $x_{k-1}$, $x_1+\cdots+x_{k-1}$, up to sign.
By Corollary~\ref{variable-symmetry}, the first $k-1$ forms on
this list all have $R^\lambda/\ideal{x_i} \cong R^\lambda/\ideal{x_1} \cong R^\nu$.
Hence $R^\nu$ can be identified intrinsically by ``majority rule'':  
it is the $\Zz$-algebra that occurs (up to isomorphism)
as the quotient of $R^\lambda$ by at least $k-1$ of the $k$ different
%%% $\pm$-equivalence classes of
primitive $k$-nilpotent linear forms
(other than the one, namely $x_1+\cdots+x_{k-1}+2x_k$, that is linearly
independent from the rest, as above).
Note that the fact that $k-1$ out of $k$ is
a well-defined ``majority'' uses the assumption that $k > 2$.

\vskip 0.1in
\noindent
{\sf Subcase IIB.} $k = 2$.

If $\lambda_2>3$, then $x_1$ is the unique primitive
$k$-nilpotent up to sign, so it is distinguished intrinsically,
as is $R^\nu \cong R/\ideal{x_1}$.

If $\lambda_2=3$, then there are two primitive $k$-nilpotents up to sign, 
namely $x_1$ and $x_1+2x_2$.
We claim that the graded $\Zz$-algebra map $\omega: R^\lambda \to R^\lambda$ defined by
    $$\omega(x_1) = x_1+2x_2, \qquad
      \omega(x_2) = -x_2, \qquad
      \omega(x_i) = x_i \ \ \text{for } 3 \leq i \leq n
    $$
is an automorphism of $R^\lambda$
interchanging $x_1$ with $x_1+2x_2$. Indeed, it is a routine calculation to check
that $\omega$ lifts to an automorphism of $\Zz[x_1,\dots,x_n]$, and
that $\omega(I_\lambda)=I_\lambda$.  In particular,
$R^\nu \cong R^\lambda/\ideal{x_1} \cong R^\lambda/\ideal{x_1+2x_2}$
may again be described up to isomorphism as the quotient of $R^\lambda$
by an arbitrary primitive $k$-nilpotent linear form.
\end{proof}

\section{The decomposable case}
\label{decomposable}
%%%%%%%%%%%%%%%%%%%%%%%%%%%%%%%%%%%%%%%%%%%%%%%%%%%%%%%%%%%%%%%%%%%%%%%%%%%%%%%%%%%%%%%

We now consider the case that $\lambda$ is decomposable, with indecomposable components
$\lambda^{(1)},\lambda^{(2)},\dots,\lambda^{(r)}$.  In this case,
$X_\lambda \cong X_{\lambda^{(1)}} \times \dots \times X_{\lambda^{(r)}}$.
Since each $X^{\lambda^{(i)}}$ has no torsion in its (co-)homology
by Theorem~\ref{Ding's-thm}, the K\"unneth formula \cite[\S 61]{Munkres} 
implies a tensor decomposition for the associated cohomology rings: 
    \begin{equation} \label{decompose-cohomology}
    H^*(X_\lambda;\:\Zz) ~\cong~ \bigotimes_{i=1}^r H^*(X_{\lambda^{(i)}};\:\Zz).
    \end{equation}
Together with the uniqueness result for indecomposable partitions
(Theorem~\ref{mainthm-indecomp}),
it would seem that we are done.  However, there is one remaining technical point:
to verify that the partitions $\lambda^{(i)}$ can be read off intrinsically 
from the structure of $H^*(X_\lambda)$ as a graded $\Zz$-algebra,
we must check that the tensor decomposition~\eqref{decompose-cohomology} is unique.

To do this, we make further use of the facts about nilpotence established
in Section~\ref{general-nilpotence-section}.  But first we must make precise
the notion of tensor decomposition, and point out how it interacts with order
of nilpotence.

%When $\lambda$ is decomposable, the Schubert variety decomposes into the products of the
%Schubert varieties corresponding to its indecomposable components $\lambda^{(i)}$.
% Similarly, the ring $R^\lambda$
%decomposes as a tensor product of the rings $R_{\lambda^i}$, as in this case, the ideal $I$
%is generated by elements
%each of which is supported only on variables from a single subpartition
%$\lambda^i$. Our goal is to uniquely
%identify this decomposition. We start by defining exactly what we mean.

For $R$ a standard graded $\Zz$-algebra, a {\it tensor decomposition} is
an isomorphism of graded $\Zz$-algebras 
$R \cong T^{(1)} \otimes \cdots \otimes T^{(r)}$ in which each $T^{(i)}$ is a standard
graded $\Zz$-algebra.  Note that any such decomposition is completely determined by
the associated direct sum decomposition of free $\Zz$-modules 
$R_1 \cong \bigoplus_{i=1}^r T^{(i)}_1$, since $T^{(i)}$ is then the
subalgebra of $R$ generated by the direct summand $T^{(i)}_1$ of $R_1$.
Say that a tensor decomposition of $R$ is
\defterm{nontrivial} if $T^{(i)} \neq \Zz$ for all $i$.  Say
$R$ is \defterm{tensor-indecomposable} if it is not $\Zz$ itself, and
has no nontrivial tensor decomposition.  

\begin{lem} \label{nilpotence-order-calc}
Suppose that $R = T^{(1)} \otimes \cdots \otimes T^{(r)}$.
Let $x \in R_1$; that is,
    $$x ~=~ x_1\otimes 1 \otimes\cdots \otimes 1 \;+\;
        1\otimes x_2 \otimes 1\cdots \otimes 1 \;+\; \dots \;+\;
        1\otimes\cdots \otimes 1 \otimes x_r$$
where $x_i \in T^{(i)}_1$.
Let $k_i$ be the nilpotence order of $x_i$.  (Recall that $k_i=1$ if and only if
$x_i=0$.)

Then the nilpotence order of $x$ is
    $$c ~=~ k_1+k_2+\dots+k_r-r+1.$$
\end{lem}

\begin{proof}
By the pigeonhole principle, each term of the multinomial expansion
of $x^c$ is divisible by $x_i^{k_i}$ for some $i$; therefore, $x^c=0$ in $R$.
For the same reason, all but one term of the multinomial expansion of $x^{c-1}$
vanishes; the exception is
    $$
    \binom{c}{k_1-1,\: \ldots,\: k_n-1} \;
    x_1^{k_1-1}\otimes x_2^{k_2-1} \otimes\cdots \otimes x_n^{k_n-1},
    $$
which is nonzero, since it is nonzero in each tensor factor.
% Here we use the fact that $R$ is a free $\Zz$-module.
\end{proof}

This calculation has immediate useful consequences.

\begin{cor} \label{nilpotent-membership}
Let $R$ be a standard graded $\Zz$-algebra with a nontrivial
tensor decomposition $R=\bigotimes_{i=1}^{r} T^{(i)}$.
Then any linear form $f \in R_1$ that achieves the minimal nilpotence 
among all elements in $R_1$ must lie in $T^{(i)}$ for some $i$.
\end{cor}

Combining Lemma~\ref{nilpotence-order-calc} with 
Proposition~\ref{variable-nilpotence-indecomp} yields the following.

\begin{cor}
\label{variable-nilpotence}
Let $\lambda$ be a partition with indecomposable components $\{\lambda^{(j)}\}_{j=1}^r$.
If $\lambda_i$ corresponds to $\lambda^{(j)}_k$ in this decomposition, then
$x_i$ is $\lambda^{(j)}_k$-nilpotent in $R^\lambda$.
\end{cor}

For example, if $\lambda$ is the decomposable partition shown in
Figure~\ref{lambda-example}, then $\lambda_1,\dots,\lambda_5$ correspond to the rows
of $\lambda^{(1)}$, and $\lambda_7,\lambda_8$ to the rows of $\lambda^{(2)}$.
Thus the variables $x_1,\dots,x_5$ have nilpotence orders $5,5,5,6,6$, respectively,
in $R^\lambda$ (and in $R^{\lambda^{(1)}}$), and $x_7,x_8$ have nilpotence orders
$2$ and $3$, respectively.  (Note that these seven variables are a
$\Zz$-basis for $R^\lambda_1$;
$x_6 \equiv -(x_1+\dots+x_5)$ does not correspond to a variable in the presentation for
$R^{\lambda^{(1)}}$.)

%With this lemma, we prove that indecomposable partitions have associated cohomology
%rings which are also indecomposable. The proof uses a result proved in
%Section~\ref{indecomposable} which states that
%the minimal nilpotence order of any element of $R^\lambda$ is $k$ for $\lambda$
%indecomposable.

\begin{prop}\label{indecomp-indecomp}
Let $\lambda$ be an indecomposable partition.  Then the ring
$R^\lambda$ is tensor-indecomposable.
\end{prop}

\begin{proof}
Let $n$ denote the number of parts in $\lambda$, and $k=\lambda_1$ its
smallest part. We proceed by induction on $n$.

If $n=1$, then clearly $R^\lambda=\Zz[x_1]/\ideal{x_1^k}$ is indecomposable.
Otherwise, suppose that $R^\lambda=T^{(1)}\otimes T^{(2)}$ is a nontrivial tensor
decomposition; we will obtain a contradiction.

By Proposition~\ref{min-nil-deg}, $x_1$ is a nilpotent of minimal order,
and hence by Corollary~\ref{nilpotent-membership}, 
without loss of generality, $x_1 \in T^{(1)}$.
Then $R^\lambda/\ideal{x_1} = T^{(1)}/\ideal{x_1} \otimes T^{(2)}$.
On the other hand, $R^\lambda/\ideal{x_1} \cong R^\nu$,
where $\nu$ is the partition obtained from 
$\lambda$ by removing the first row and column.
Since $\lambda$ is indecomposable, so is $\nu$.
By the inductive hypothesis, the decomposition 
$T^{(1)}/\ideal{x_1} \otimes T^{(2)}$ must be trivial; 
that is, $T^{(1)}/\ideal{x_1} \cong \Zz$,
and $T^{(1)}$ must be generated
by $x_1$ as a $\Zz$-algebra, i.e., $T^{(1)}=\Zz[x_1]/\ideal{x_1^k}$.
Therefore, exactly one member of the set
    $$L = \{x_2+\alpha x_1 : \alpha \in \Zz\}$$
belongs to $T^{(2)}_1$.  Let $\ell$ be the nilpotence order of that one form;
then all other elements of $L$ have nilpotence order
$k+\ell-1>\ell$ by Lemma~\ref{nilpotence-order-calc}.
Let $m=\lambda_2$; note that $m\geq 3$ since $\lambda$ is indecomposable.
By Proposition~\ref{Leray-Hirsch} we can work in the algebra 
$R^{(\lambda_1,\lambda_2)} = R^{(k,m)}$, 
namely the quotient of $\Zz[x_1,x_2]$ by the ideal
    $$\ideal{G_1 = x_1^k,\qquad
    G_2 = x_2^{m-1} + x_2^{m-2} x_1 + \ldots + x_2^{m-k} x_1^{k-1}}.$$

Let $\alpha \in \Zz$ be arbitrary.  We will show that no linear form
$x_2+\alpha x_1$ has nilpotence order strictly less than $m$.  Indeed,
  \begin{eqnarray*}
    (x_2+\alpha x_1)^{m-1}
    &=& \sum_{j=0}^{m-1} \binom{m-1}{j} x_2^j \alpha^{m-j-1} x_1^{m-j-1}\\
    &=& \left(\sum_{j=0}^{m-2} \binom{m-1}{j} x_2^j \alpha^{m-j-1} x_1^{m-j-1}\right) + G_2
        - \sum_{j=0}^{m-2} x_2^j x_1^{m-j-1}\\
    &=& G_2 + \sum_{j=0}^{m-2} \left(\binom{m-1}{j}\alpha^{m-j-1}-1\right)
        x_2^j x_1^{m-j-1}\\
    &\equiv& \sum_{j=m-k}^{m-2} \left(\binom{m-1}{j}\alpha^{m-j-1}-1\right)
        x_2^j x_1^{m-j-1} \quad \mod{\ideal{G_1, G_2}}.
  \end{eqnarray*}
This last expression is exactly the standard form of
$(x_2+\alpha x_1)^{m-1}$.  For $j=m-2$,
the summand is $((m-1)\alpha-1) x_2^{m-2} x_1$; since $m>3$ and $\alpha$
is an integer, the coefficient is nonzero.
Therefore $(x_2+\alpha x_1)^{m-1} \neq 0$.

On the other hand, $x_2^m = 0$ in $R^\lambda$ by 
Proposition~\ref{variable-nilpotence-indecomp}.
Therefore $x_2$ must be the unique element of $L$ with minimal nilpotence order
$m=\ell$, and every other element of $L$ must have nilpotence order $k+m-1$.
But there are no standard monomials in $x_1,x_2$ of degree greater than $(k-1)+(m-2)
= k+m-3$, which implies that every element of $L$ has nilpotence
order $k+m-2$ or less.  This contradiction completes the proof.
\end{proof}

We now establish the key fact of the decomposable case, that these
decompositions are actually unique.

\begin{lem}\label{keylemma}
The ring $R^\lambda$ has a unique tensor decomposition
into tensor-indecomposables.  Specifically, if
$\lambda$ has indecomposable components
$\lambda^{(1)},\: \lambda^{(2)},\: \ldots,\: \lambda^{(r)}$,
then
$$
R^\lambda = R^{\lambda^{(1)}} \otimes \cdots \otimes R^{\lambda^{(r)}}
$$
is the unique tensor decomposition of $R^\lambda$, up to permuting the factors.
\end{lem}

\begin{proof}
The existence is immediate, since each
$R^{\lambda^{(i)}}$ is tensor-indecomposable by Lemma~\ref{indecomp-indecomp}.
For uniqueness, we proceed by induction on the number of rows of $\lambda$.
If $\lambda$ has only one row, the statement is trivial.

Suppose that $R^\lambda = \otimes_{i=1}^s T^{(i)}$
is a tensor decomposition with each $T^{(i)}$ tensor-indecomposable, so that
  \begin{subequations}
  \begin{align}
  \bigotimes_{i=1}^s T^{(i)} &\quad=\quad R^\lambda
    \quad=\quad \bigotimes_{j=1}^r R^{\lambda^{(j)}},
    \label{factorization-tensor} \\
  \bigoplus_{i=1}^s T^{(i)}_1 &\quad=\quad R^\lambda_1 
    \quad=\quad \bigoplus_{j=1}^r R^{\lambda^{(j)}}_1.
    \label{factorization-directsum}
  \end{align}
  \end{subequations}
Let $k$ be the minimal nilpotence order of any element of $R^\lambda_1$.
Then $k = \min\{\lambda^{(j)}_1: 1 \leq j \leq r\}$
by Corollary~\ref{variable-nilpotence}.
Without loss of generality, we may re-index so that
$k=\lambda^{(1)}_1$; then $x_1$ is
a linear form of nilpotence order $k$.
By Corollary~\ref{nilpotent-membership},
$x_1$ must belong to one of the $T^{(i)}$, say $T^{(1)}$.
Let $\nu, \nu^{(1)}$ be the partitions obtained by removing the left column and bottom
row of $\lambda, \lambda^{(1)}$, respectively.
Then
  \begin{subequations}
  \begin{align}
  T^{(1)}/\langle x_1 \rangle \otimes \bigotimes_{i=2}^s T^{(i)} 
    &\quad=\quad R^\lambda/\langle x_1 \rangle \quad=\quad
    R^{\nu^{(1)}} \otimes \bigotimes_{j=2}^r R^{\lambda^{(j)}},
    \label{quotient-factorization-tensor} \\
  T^{(1)}_1/\Zz x_1 \oplus \bigoplus_{i=2}^s T^{(i)}_1 
    &\quad=\quad R^\lambda_1/\Zz x_1 \quad=\quad
    R^{\nu^{(1)}}_1 \oplus \bigoplus_{j=2}^r R^{\lambda^{(j)}}_1. 
    \label{quotient-factorization-directsum}
  \end{align}
  \end{subequations}
By induction, the rightmost expression in~\eqref{quotient-factorization-tensor}
is the unique tensor decomposition of $R^\nu$ into tensor-indecomposables
(possibly with a superfluous factor $R^{\nu^{(1)}} \!\! = \Zz$
if $\lambda^{(1)}$ has only one part).  Thus the rightmost expression
in~\eqref{quotient-factorization-directsum} is unique---clearly
not as a direct sum decomposition
of $R^\lambda_1/\Zz x_1$ as a $\Zz$-module, but as a direct sum decomposition
which induces a tensor decomposition of $R^\lambda/\ideal{x_1}$.

Now assume that $\lambda^{(1)}$ has $m$ rows, so that
$x_1,x_2,\ldots,x_m$ generate $R^{\lambda^{(1)}}$
as a $\Zz$-subalgebra of $R^\lambda$.  For each $\ell \in \{2,\dots,m\}$,
consider the image $\bar x_\ell$ of $x_\ell$ in $R^\nu_1 =R^\lambda_1/\Zz x_1$.
Since each $\bar x_\ell$ belongs to the direct summand
$R^{\nu^{(1)}}_1$ on the left side of the unique
decomposition~\eqref{quotient-factorization-directsum}, it must belong either
to $T^{(1)}_1/\Zz x_1$, or to $T^{(i)}_1$ for
some $i \geq 2$.  On the other hand, Corollary~\ref{variable-nilpotence}
tells us that $x_\ell$ is $\lambda^{(1)}_\ell$-nilpotent in $R^\lambda$,
but $\bar x_\ell$
is $\nu^{(1)}_{\ell-1}$-nilpotent in $R^\nu$.  That is, the nilpotence 
order 
of $x_\ell$
drops by $1$ in the quotient by $x_1$ (because $\nu^{(1)}_{\ell-1} = 
\lambda^{(1)}_{\ell}-1$).
If $\bar x_\ell \in T^{(i)}_1$ for some $i \geq 2$, then this
last observation contradicts Lemma~\ref{nilpotence-order-calc}.  
Therefore $\bar x_\ell \in T^{(1)}_1/\Zz x_1$, from which we
conclude that $T^{(1)}_1/\Zz x_1 \supseteq R^{\nu^{(1)}}_1$.

Consequently, the uniqueness property of the
decomposition~\eqref{quotient-factorization-directsum} implies that
    $$
    T^{(1)}_1/\Zz x_1 
      ~=~ R^{\nu^{(1)}}_1 \oplus \bigoplus_{u \in U} R^{\lambda^{(j_u)}}_1
    $$
for some subset $U \subset \{2,3,\ldots,r\}$.  
Since $x_1$ lies in both $T^{(1)}$ and $R^{\lambda^{(1)}}$, we conclude that
    $$
    T^{(1)}_1 ~=~ R^{\lambda^{(1)}}_1 \oplus \bigoplus_{u \in U} R^{\lambda^{(j_u)}}_1
    $$
and, since $T^{(1)}$ is a standard graded $\Zz$-algebra,
    $$
    T^{(1)} ~=~ R^{\lambda^{(1)}} \otimes \bigotimes_{u \in U} R^{\lambda^{(j_u)}}.
    $$
But $T^{(1)}$ was assumed to be indecomposable, so this forces $U =\emptyset$.
Hence $T^{(1)}_1 = R^{\lambda^{(1)}}_1$ and
$T^{(1)}_1/\Zz x_1 = R^{\nu^{(1)}}_1$.  By the uniqueness property
of~\eqref{quotient-factorization-directsum}, we must have $r=s$, and after re-indexing,
$T^{(i)}_1 = R^{\lambda^{(i)}}_1$ for $i=2,3,\ldots,r$.  Thus
the two tensor decompositions in \eqref{factorization-tensor} are identical.
\end{proof}

The nontrivial implication (iii) $\implies$ (i) in
the main result, Theorem~\ref{main-theorem}, is now immediate from
Lemma~\ref{keylemma} and Theorem~\ref{mainthm-indecomp}.

%We close with an example that illustrates why it was important to use
%integer coefficients in the cohomology rings $H^*(X_\lambda)=R^\lambda$,
%and why some of the technical work of Section~\ref{general-nilpotence-section}
%was required.

\begin{remark}
\label{invertible-2}
As we shall now demonstrate, it was essential to
study the cohomology of $X_\lambda$ with integer coefficients.
If $A$ is a coefficient ring in which
$2$ is invertible, then
Proposition~\ref{indecomp-indecomp}, Lemma~\ref{keylemma}
and Theorem~\ref{main-theorem} would all fail to hold
if ``graded $\Zz$-algebra'' was replaced with ``graded $A$-algebras''.
That is, Ding's Schubert varieties are \emph{not}
classified up to isomorphism by their cohomology with $A$-coefficients.
For example, consider the indecomposable partition $\lambda=(2,3)$.
By completing the square, one has
    \begin{align*}
    R^{(2,3)} \otimes_\Zz A
    ~&\cong~ A[x_1,x_2]\:/\ideal{x_1^2,\, x_2^2+x_1 
x_2+\frac{1}{4}x_1^2}\\
    ~&=~     A[x_1,x_2]\:/\ideal{x_1^2,\, \left( x_2+\frac{1}{2} x_1 \right)^2}\\
    ~&\cong~ A[x_1]\:/\ideal{x_1^2} \,\otimes\, A[y]\:/\ideal{y^2}.
    \end{align*}
Thus indecomposable partitions do not lead to tensor-indecomposable
graded $A$-algebras.  This also leads to ``extra'' isomorphisms
among the cohomology rings $H^*(X_\lambda;A) \cong R^\lambda \otimes_\Zz A$.
For example, the partition $\mu=(2,2,4)$ has indecomposable components
$\mu^{(1)}=\mu^{(2)}=(2)$.  Since $R^{(2)} \cong \Zz[x]/\ideal{x^2}$, one has
    $$
    R^\mu \otimes_\Zz A ~\cong~ A[x]\:/\ideal{x^2} \otimes_A A[x]\:/\ideal{x^2} 
    ~\cong~ R^\lambda \otimes_\Zz A
    $$ 
even though $\lambda=(2,3)$ and $\mu=(2,2,4)$ do not have the
same indecomposable partition components.
\end{remark}

%%%%%%%%%%%%%%%%%%%%%%%%%%%%%%%%%%%%%%%%%%%%%%%%%%%%%%%%%%%%%%%%%%%%%%%%%%%%%%%%%%%%%%%

\end{document}